\nonstopmode
\documentclass{m2an}
\usepackage[T1]{fontenc}
\usepackage{lmodern}
\usepackage[utf8]{inputenc}
\usepackage{array}
\usepackage{textcomp}
\usepackage{mathrsfs}
\usepackage{amssymb}
\usepackage{amsmath}
\usepackage{graphicx}
\usepackage{cite}
\usepackage[colorlinks=true]{hyperref}
\DeclareMathOperator*{\argmax}{argmax}
\DeclareMathOperator*{\arginf}{arginf}
\DeclareMathOperator*{\esssup}{esssup}
\newcommand{\N}{{N}}

\newcommand{\R}{\mathbf{R}}

\newcommand{\ud}{\mathrm{d}}
\newcommand{\delt}{\Delta t}
\newcommand{\utk}{{\widetilde u}^k}
\newcommand{\ut}[1]{{\widetilde u}^{#1}}
\newcommand{\rhog}{\text{\boldmath{$\rho$}}}
\newcommand{\intzu}{\int_0^1}

\newcommand{\phiuz}{\Phi^{u_0}_l(x)}

\newcommand{\phifsp}{\Phi^{fS}_p(\cdot)}
\newcommand{\phiuzp}{\Phi^{u_0}_l(\cdot)}

\newcommand{\phiftr}{\Phi^{fT}_l}
\newcommand{\phiftrt}{\Phi^{fT}_l}

\renewcommand{\epsilon}{\varepsilon}

\newcommand{\partieneg}[1]{\left[ #1 \right]_-}
\everymath{\displaystyle\everymath{}}
\providecommand{\abs}[1]{\left\lvert#1\right\rvert}
\providecommand{\norm}[1]{\left\lVert#1\right\rVert}
\usepackage{float}
\floatstyle{boxed}
\newfloat{algo}{h}{loa}
\floatname{algo}{Algorithm}
\numberwithin{equation}{section}
\usepackage{ifthen}
\newboolean{light}
\lighttrue
\begin{document}
\title{Certified reduced-basis solutions of viscous Burgers equation parametrized by initial and boundary values}
\begin{abstract} We present a reduced basis offline/online procedure for viscous Burgers initial boundary value problem, enabling efficient approximate computation of the solutions of this equation for parametrized viscosity and initial and boundary value data. This procedure comes with a fast-evaluated rigorous error bound certifying the approximation procedure. Our numerical experiments show significant computational savings, as well as efficiency of the error bound. \end{abstract}
\author{Alexandre Janon} \address{Joseph Fourier University, LJK/MOISE, BP 53, 38041 Grenoble Cedex, France ; \email{alexandre.janon@imag.fr\ \& maelle.nodet@inria.fr\ \& clementine.prieur@imag.fr}. }
\author{Maëlle Nodet} \sameaddress{1}
\author{Clémentine Prieur} \sameaddress{1}
\maketitle
\section*{Introduction}
This paper is set in the context of sensitivity analysis and uncertainty analysis in geophysical models. Such models typically involve a wide range of parameters, such as:
source terms (climatic forcings, heat/wind/matter fluxes),
 boundary conditions (forcings, open boundaries),
and the initial state of the system.

Their study generally leads to parametrized partial differential equations (PDEs). These equations often involve poorly-known parameters. Therefore, it is important to be able to measure the impact of a given parameter on the quality of the solution, and also to identify the "sensitive" parameters, that is, the parameters for which a small variation implies a large variation of the model solution. Due to their ability to perform global sensitivity analyses for nonlinear models, stochastic tools \cite{helton2006survey,saltelli-sensitivity} are rapidly expanding. These methods require "many queries," that is solving the parametrized PDE for a large (say, thousands) number of values of the parameters. When analytic solution to the PDE is not known (as it is often the case), one has to use a numerical method (such as finite difference or finite element) to compute an approximate value of the solution. Such methods lead to computer codes that could take a large time to produce an accurate-enough approximation --- for a single value of the parameter. Having the "many-query" problem in mind, it is crucial to design a procedure that solves the equation for several values of the parameter faster than the naïve approach of calling the numerical code for each required instance of the parameter. 

The reduced basis (RB) method is such a procedure; we split the overall computation into two  successive parts: one part, the \emph{offline} phase, makes use of the standard, computationally intensive numerical procedure used to solve the PDE to gather "knowledge" about solutions of the latter; and the other one, the \emph{online} phase, where we rely on data collected during the offline phase to compute, for each desired instance of the parameter, a good approximation of the solution, for a per-instance cost that is orders of magnitude smaller than the cost of one run of the standard numerical code. The advantage is that, for a sufficiently large number of online evaluations, the fixed cost of the offline phase will be strongly dominated by the reduction in the marginal cost provided by the online procedure. This cost reduction is made possible by the fact that, in most cases, the desired solutions of the PDE, for all the considered values of the parameter, lie in some manifold of functions that is "close" to a low-dimensional linear subspace. One goal of the offline phase is to find such a suitable subspace, so that the online procedure can look for the solution of the PDE as an element of the subspace --- so as to reduce the number of degrees of freedom and thus the computational cost. One interesting feature of the RB approach is that it comes with an \emph{online error bound}, that is a (provably) certified, natural norm, fast-computed (i.e. almost of the same complexity of the online phase) upper bound of the distance between the solution provided by the online phase (called the \emph{reduced}, or \emph{online}) solution and the one given by the standard, expensive numerical procedure (called the \emph{full} or \emph{reference} solution). This "certified RB" framework has been developed for \emph{affinely} pa\-ra\-me\-tri\-zed second-order elliptic linear PDEs in \cite{nguyen2005certified}. It has been extended to nonlinear, non-affinely parametrized, parabolic PDEs, see e.g. \cite{grepl2005posteriori}, \cite{grepl2007efficient} and applied to problems such as steady incompressible Navier-Stokes \cite{veroy2005certified}. Moreover theoretical work has been done to ensure \emph{a priori} convergence of the RB procedure \cite{buffa2009apriori}.

In this paper, we are interested in the RB reduction of the time-dependent viscous Burgers equation (which will serve as a "test case" for the "real" equations modelling geophysical fluids we are interested in). Papers \cite{haasdonk2008reduced} and \cite{rovas2006reduced} extend certified RB methodology to linear initial-boundary value problems. The case of homogeneous Dirichlet boundary conditions, zero initial value and fixed (\emph{i.e.}, not parametrized) source term has been treated in \cite{veroy2003reduced}, \cite{nguyen2009reduced} ; in these works, the only parameter was the viscosity coefficient. Parametrization of initial and Dirichlet boundary conditions (treated using a conversion to a homogeneous Dirichlet problem) has been done in \cite{jung2009reduced}, for a general multidimensional quadratically nonlinear equation. Our methodology allows parametrization of the viscosity, the initial state, the source term and and of the boundary conditions. Compared to the works cited above, we use a weak (penalization) treatment of the Dirichlet boundary conditions. We will see that this weak treatment is more favorable in terms of both computation and storage complexities. Besides, our paper features a new error bound, which has shown to be, in all the testcases we performed, much more efficient than the existing bound. In particular, we will show that our bound exhibits improved sharpness for low viscosities, where \cite{nguyen2009reduced} pointed out the moderate efficiency of the presented \emph{a posteriori} bound. 

This paper is organised as follows: in the first part, we introduce the viscous Burgers equation, and present a standard numerical procedure used to solve it; in the second part, we expose our offline/online reduction procedure; in the third part, we develop a certified online error bound; finally in the fourth part we validate and discuss our results based on numerical experiments. 

\section{Model}
In this section, we describe the model we are interested in. Subsection \ref{ss:eqn} introduces the viscous Burgers equation, while Subsection \ref{ss:num} presents the "full" numerical procedure on which our reduction procedure, described in Section \ref{s:redproc}, relies on.
\subsection{Equation}
\label{ss:eqn}
We are interested in $u$, function of space $x\in[0;1]$ and time $t\in[0;T]$ (for $T>0$), with regularity: \\ \mbox{$u \in C^1 \left([0,T], H^1(]0,1[) \right)$},
satisfying the \emph{viscous Burgers equation}:
\begin{equation}
\label{e:burg}
\frac{\partial u}{\partial t} + \frac{1}{2} \frac{\partial}{\partial x}(u^2) - \nu \frac{\partial^2 u}{\partial x^2} = f 
\end{equation}
where $\nu\in\R^+_*$ ($\R$ denotes the set of real numbers, $\R^+_*$ the set of positive real numbers) is the \emph{viscosity}, and $f \in C^0\left([0,T], L^2(]0,1[)\right)$ is the \emph{source term}.

For $u$ to be well-defined, we also prescribe initial values $u_0 \in H^1(]0,1[)$:
\begin{equation}
\label{e:init}
u(t=0,x)=u_0(x) \;\; \forall x\in[0;1]
\end{equation}
and boundary values $b_0, b_1 \in C^0([0,T])$:
\begin{equation}
\label{e:boundd}
\begin{cases} u(t,x=0)=b_0(t) \\ u(t,x=1)=b_1(t) \end{cases} \;\; \forall t\in[0;T]  
\end{equation}
Where $b_0$, $b_1$ and $u_0$ are given functions, supposed to satisfy \emph{compatibility conditions}:
\begin{equation}
\label{e:compatcond}
u_0(0)=b_0(0) \;\; \text{ and } \;\; u_0(1)=b_1(0) 
\end{equation}

This problem can be analyzed by means of the Cole-Hopf substitution (see \cite{hopf1950partial} for instance), which turns \eqref{e:burg} into the heat equation, leading to an integral representation of $u$.

\subsection{Numerical resolution}
\label{ss:num}
We now describe the "expensive" numerical resolution of the problem described above that will serve as our reference for the reduction procedure described in the next section. We proceed in two steps: space discretization in paragraph \ref{sss:space} and time discretization in paragraph \ref{ss:numT}.
\subsubsection{Space discretization}
\label{sss:space}
For space discretization, we use a $\mathbf{P}^1$ finite element procedure with weak (penalty) setting of the Dirichlet boundary conditions \eqref{e:boundd}.

We first have to write the weak formulation of our PDE ; to do so, we multiply \eqref{e:burg} by a function $v\in H^1(]0;1[)$ and integrate over $]0;1[$:
\begin{multline}
\label{e:wf1}
\intzu \frac{\partial u}{\partial t}(t,x) v(x) \ud x + \frac{1}{2} \intzu \frac{\partial (u^2)}{\partial x}(t,x) v(x) \ud x \\
- \nu \intzu \frac{\partial^2 u}{\partial x^2}(t,x) v(x) \ud x = \intzu f(t,x) v(x) \ud x \;\; \forall v \in H^1(]0;1[) \; \forall t \in [0;T] 
\end{multline}
Next, we integrate by parts the second and the third integral appearing in the left hand side of the previous equation:
\begin{align*}
\intzu \frac{\partial (u^2)}{\partial x}(t,x) v(x) \ud x &= - \intzu u^2(t,x) \frac{\partial v}{\partial x}(x) \ud x + \left[ u^2(t,\cdot) v \right]_0^1 \\
\intzu \frac{\partial^2 u}{\partial x^2}(t,x) v(x) \ud x &= - \intzu \frac{\partial u}{\partial x}(t,x) \frac{\partial v}{\partial x}(x) \ud x + \left[ \frac{\partial u}{\partial x}(t,\cdot) v \right]_0^1
\end{align*}
Inserting this into \eqref{e:wf1}, we get:
\begin{multline}
\label{e:wf2}
\intzu \frac{\partial u}{\partial t}(t,x) v(x) \ud x - \frac{1}{2} \intzu u^2(t,x) \frac{\partial v}{\partial x}(x) \ud x
+ \nu \intzu \frac{\partial u}{\partial x}(t,x) \frac{\partial v}{\partial x}(x) \ud x \\
+ \frac{1}{2} \left[ u^2(t,\cdot) v \right]_0^1
- \nu \left[ \frac{\partial u}{\partial x}(t,\cdot) v \right]_0^1 
=
\intzu f(t,x) v(x) \ud x \;\; \forall v \in H^1(]0;1[) \; \forall t \in [0;T] 
\end{multline}
To get rid of the two boundary terms arising in the integrations by parts, one usually restricts $v$ to satisfy $v(0)=v(1)=0$ so as to make the boundary terms disappear; the Dirichlet boundary conditions \eqref{e:boundd} are then incorporated "outside" of the weak formulation. However, the reduction framework we are to expose later requires the boundary conditions to be ensured by the weak formulation itself.  The Dirichlet penalty method, presented in \cite{babuska1973finite}, is a way of doing so, at the expense of a slight approximation error. This method entails replacement of boundary conditions \eqref{e:boundd} with the following conditions:
\begin{equation}
\label{e:bound2}
\begin{cases}
- \frac{1}{2} u^2(t,x=0)+\nu \frac{\partial u}{\partial x} (t,x=0) = P \left( u(t,x=0)-b_0(t) \right) \\ 
\frac{1}{2} u^2(t,x=1)-\nu \frac{\partial u}{\partial x} (t,x=1) = P \left( u(t,x=1)-b_1(t) \right) 
\end{cases} \;\; \forall t\in[0;T]
\end{equation}
with a fixed \emph{penalization constant} $P>0$.

The intuitive idea underlying \eqref{e:bound2} is that it can clearly be rewritten as:
\begin{equation*}
\begin{cases}
u(t,x=0) = b_0(t) + \frac{1}{P} \left( - \frac{1}{2} u^2(t,x=0)+\nu \frac{\partial u}{\partial x} (t,x=0) \right) \\
u(t,x=1) = b_1(t) + \frac{1}{P} \left( \frac{1}{2} u^2(t,x=1)-\nu \frac{\partial u}{\partial x} (t,x=1) \right) 
\end{cases} \;\; \forall t\in[0;T]
\end{equation*}
so that \eqref{e:boundd} is asymptotically verified for $P \rightarrow +\infty$. The reader can refer to \cite{barrett1986finite} for rigorous \emph{a priori} error estimates when using Dirichlet penalty in the linear elliptic case. 

In practice, we can check if our approximation is sufficiently accurate by means of the following \emph{a posteriori} procedure: we take for $P$ some large value (typically $P=10^7$), we compute (numerically) an approximate solution $u_d$, using the procedure we are currently describing, and we check if an indicator of the amount of failure in verification of \eqref{e:boundd} is small enough; such an indicator can be, for instance:
\begin{equation}
\label{e:posterioridiagnosis}
\epsilon_b = \sup_t \left[ \max \left( | u_d(t,x=0)-b_0(t) |, |u_d(t,x=1)-b_1(t)| \right) \right] 
\end{equation}
where the supremum is taken over all discrete time steps. If this indicator is larger than a prescribed tolerance, then $P$ has to be increased. Our numerical results in Section \ref{s:numres} will assert this condition. We can then invoke the well-posedness of the boundary/initial value problem (specifically, continuous dependence on the boundary values) to ensure that the solution of \eqref{e:burg}, \eqref{e:init} and \eqref{e:bound2} will be close to the solution of \eqref{e:burg}, \eqref{e:init} and \eqref{e:boundd}. This reasoning is analogous to the one made when omitting the approximation made when replacing exact boundary values by their discretized counterparts. 

Going back to our weak formulation, we multiply the first line of \eqref{e:bound2} by $v(0)$, the second one by $v(1)$ and add up these two equations. We get that:
\begin{multline*}
\frac{1}{2} \left[ u^2(t,\cdot) v \right]_0^1 - \nu \left[\frac{\partial u}{\partial x}(t,\cdot) v\right]_0^1
=
P \left[ \left( u(t,x=0) v(0) - b_0(t) v(0) \right) + \left( u(t,x=1) v(1) - b_1(t) v(1) \right) \right]
\end{multline*}
Putting it back into \eqref{e:wf2} and isolating the terms not involving $u$ on the right-hand side yields the following weak formulation:
\begin{multline}
\label{e:wf3}
\intzu \frac{\partial u}{\partial t} v - \frac{1}{2} \intzu u^2 \frac{\partial v}{\partial x} + \nu \intzu \frac{\partial u}{\partial x} \frac{\partial v}{\partial x} + P \left( u(t,x=0)v(0) + u(t,x=1)v(1) \right)
\\ = \intzu f(t,\cdot) v + P \left( b_0(t) v(0) + b_1(t) v(1) \right) \;\; \forall v \in H^1(]0;1[) \; \forall t \in [0;T] 
\end{multline}
that is:
\begin{multline}
\label{e:wf4}
\Big\langle \frac{\partial u}{\partial t}(t,\cdot),v \Big\rangle  + c(u(t,\cdot),u(t,\cdot),v) + \nu a (u(t,\cdot),v) + B(u(t,\cdot),v) \\
 = \ell(v,t) + b_0(t) \beta_0(v) + b_1(t) \beta_1(v) \;\; \forall v \in H^1(]0;1[), \forall t \in [0,T]
\end{multline}
by introducing the following notations (for all $v,w,z\in H^1(]0;1[)$, $t\in[0;T]$):
\begin{align*}
	\langle w,v \rangle &= \intzu wv & 
a(w,v) &= \intzu \frac{\partial w}{\partial x} \frac{\partial v}{\partial x} \\
B(w,v) &= P \left( w(0) v(0) + w(1) v(1) \right) &
\ell(v,t) &= \intzu f(t,\cdot) v \\
\beta_0(v) &= P v(0) &
\beta_1(v) &= P v(1) 
\end{align*}
and:
\[ c(w,v,z) = - \frac{1}{2} \intzu wv \frac{\partial z}{\partial x} \]
for every $w$, $v$ and $z$ in $H^1(]0;1[)$ for which $wv\frac{\partial z}{\partial x} \in L^1(]0;1[)$.

The weak formulation is then discretized with Lagrange $\mathbf{P}^1$ finite elements (see \cite{quarteroni2008numerical}) by choosing some integer $\mathcal N$ and considering a uniform subdivision of $[0;1]$ with $\mathcal N+1$ nodes: $\left\{ x_i 
\right\}_{i=0,\ldots,\mathcal N}$ and, for each $i=0,\ldots,\mathcal N$ we denote by $\phi_i$ the piecewise-affine ``hat'' function whose value is $1$ on $x_i$ and $0$ on every other nodes. 

We denote by $X$ the linear subspace of $H^1(]0;1[)$ spanned by $\left\{ \phi_i \right\}_{i=0,\ldots,\mathcal N}$. We also set $\norm{\cdot}$ to be the $L^2(]0;1[)$ norm.

Every $\psi \in X$ can be written as $\psi = \sum_{j=0}^{\mathcal N} \psi_j \phi_j$, and we have $\psi(x_i)=\psi_i$ for all $i=0,\ldots,\mathcal N$. This justifies that $\pi$ defined below is a projection of $H^1(]0;1[)$ onto $X$:
\[ 
\pi : \left\{ \begin{aligned} H^1(]0;1[) &\rightarrow X \\ \psi &\mapsto \sum_{j=0}^{\mathcal N} \psi(x_j) \phi_j \end{aligned} \right.
\]

The space discretization of our problem is the following: for all $t\in [0;T]$, find $u(t,\cdot) \in X$ so that :
\begin{equation}
\label{e:spdisc}
\begin{cases}
u(t=0,\cdot) = \pi(u_0) \\
\Big\langle \frac{\partial u}{\partial t}(t,\cdot),v \Big\rangle  + c(u(t,\cdot),u(t,\cdot),v) + \nu a (u(t,\cdot),v) + B(u(t,\cdot),v) \\
\hspace*{4cm} =  \ell_{\pi}(v,t) + b_0(t) \beta_0(v) + b_1(t) \beta_1(v) \;\; \forall v \in X 
\end{cases}
\end{equation}
where 
\[ \ell_{\pi}(v,t) = \intzu \pi\left(f(t,\cdot)\right) v \]
Note that $u$ now stands for a discrete solution, while it was used to designate an analytic solution before. 

\subsubsection{Time discretization}
\label{ss:numT}
We now discretize \eqref{e:spdisc} in time using the backward Euler scheme: we choose a timestep $\delt>0$ and consider an uniform subdivision of $[0;T]$: $\left\{ t_k = k \delt \right\}_{k=0,\ldots,\mathcal T}$ where $\mathcal T = \frac{T}{\delt}$.

Our fully discrete problem is: for $k=0,\ldots,\mathcal T$, find $u^k \in X$, approximation of $u(t_k,\cdot)$, satisfying:
\begin{subequations}
\label{e:disc}
\begin{equation}
\label{e:disc1}
u^0 = \pi(u_0) 
\end{equation}
\text{and: }
\begin{multline}
\label{e:disc2}
\Big\langle \frac{u^k-u^{k-1}}{\delt},v \Big\rangle  + c(u^k,u^k,v) + \nu a (u^k,v) + B(u^k,v) \\
\hspace*{2cm} =  \ell_\pi(v,t_k) + b_0(t_k) \beta_0(v) + b_1(t_k) \beta_1(v) \;\; \forall v \in X \; \forall k=1,\ldots,\mathcal T
\end{multline}
\end{subequations}
We sequentially compute $\{u^k\}_{k=0,\ldots,\mathcal T}$ in the following way: $u^0$ comes straightforwardly from  \eqref{e:disc1}, and for $k=1,\ldots,\mathcal T$, $u^k$ depends on $u^{k-1}$ through \eqref{e:disc2}, which can be rewritten:
\begin{multline}
\label{e:nleq}
\frac{1}{\delt} \langle u^k,v \rangle + c(u^k,u^k,v) + \nu a (u^k,v) + B(u^k,v) \\ = \frac{1}{\delt} \langle u^{k-1},v \rangle + \ell_\pi(v,t_k) + b_0(t_k) \beta_0(v) + b_1(t_k) \beta_1(v) \;\; \forall v \in X 
\end{multline}
We can now expand our unknown $u^k \in X$ on the $\{\phi_j\}_j$ basis:
$ u^k = \sum_{j=1}^{\mathcal N} u^k_j \phi_j, $
and the vector $\left( u^k_j \right)_j$ becomes our new unknown.

Moreover, it is sufficient for (linear-in-$v$) relation \eqref{e:nleq} to be satisfied for all $v$ in a basis of $X$, namely for $v=\phi_i, \; \forall i=0,\ldots,\mathcal N$.

So \eqref{e:nleq} can be rewritten as a nonlinear (due to the nonlinearity in $c(u^k,u^k,v)$) system of $\mathcal N + 1$ equations (one for each instantation $v=\phi_i$) involving $\left( u^k_j \right)_{j=0,\ldots,\mathcal N}$. This nonlinear system is solved using Newton iterations: starting with an initial guess $\overline{u^k}$, one looks for $\delta = \sum_{j=1}^{\mathcal N} \delta_j \phi_j$ so that $u^k = \overline{u^k} + \delta$ satisfies the linearization near $\delta=0$ of \eqref{e:nleq} for $v=\phi_i$, $i=0,\ldots,\mathcal N$, that is to say:
\begin{multline}
\label{e:lzedeq}
\frac{1}{\delt} \langle \overline{u^k}+\delta,\phi_i \rangle + c(\overline{u^k},\overline{u^k},\phi_i) + 2 c(\overline{u^k},\delta,\phi_i) + \nu a (\overline{u^k}+\delta,\phi_i) + B(\overline{u^k}+\delta,\phi_i) \\ = \frac{1}{\delt} \langle u^{k-1},\phi_i \rangle + \ell_\pi(\phi_i,t_k) + b_0(t_k) \beta_0(\phi_i) + b_1(t_k) \beta_1(\phi_i) \;\; \forall i=0,\ldots,\mathcal N 
\end{multline}
System \eqref{e:lzedeq} is a $(\mathcal N+1) \times (\mathcal N+1)$ linear system involving $\left( \delta_j \right)_{j=0,\mathcal N}$. Once solved for $\delta$, one can test for convergence of the Newton iteration: if the norm of $\delta$ is smaller than a prescribed precision, then we stop here, produce $u^k$ and get to the next time step $k+1$; otherwise we do one more Newton step, this time using $u^k$ as initial guess $\overline{u^k}$.

We can note that the linear system to be solved at each Newton step is not symmetric but is sparse. Due to this sparsity, its solution (using an iterative method such as BICGSTAB or GMRES) takes about $O\left( \mathcal N^2 \right)$ operations in the worst case. One can take advantage of the tridiagonal structure of the matrix (which is present in the one-dimensional case, since $\phi_i$ and $\phi_j$ have no common support if $|j-i|>1$) and use Thomas' algorithm \cite{strikwerda2004finite} to solve the system with $O(\mathcal N)$ operations.

\section{Reduction procedure}
\label{s:redproc}
In this section, we show the offline/online procedure announced in the introduction to produce reduced basis solutions of the problem formed by \eqref{e:burg}, \eqref{e:init} and \eqref{e:boundd}, based on the "full basis" numerical method presented above. We begin by a description of our parameters in Subsection \ref{ss:params}. Our offline/online reduction procedure is described subsequently, in Subsection \ref{ss:redbcalc}.

\subsection{Parameters}
\label{ss:params}
We parametrize $u_0$, $b_0$, $b_1$ and $f$ as:
\begin{align*}
	b_0(t) &= b_{0m} + \sum_{l=1}^{n(b_0)} A^{b_0}_l \Phi^{b_0}_l(t) & 
	b_1(t) &= b_{1m} + \sum_{l=1}^{n(b_1)} A^{b_1}_l \Phi^{b_1}_l(t) \\
	f(t,x) &= f_m + \sum_{l=1}^{n_T(f)} \sum_{p=1}^{n_S(f)} A^f_{lp} \Phi^{fT}_l(t) \Phi^{fS}_p(x) &
	u_0(x) &= u_{0m} + \sum_{l=1}^{n(u_0)} A^{u_0}_l \Phi^{u_0}_l(x)
\end{align*}
The functions involved in these linear combinations $\Phi^{b_0}$, $\Phi^{b_1}_l$, $\Phi^{fT}_l$, $\Phi^{fS}_p$ and $\Phi^{u_0}_l$ as well as the number of terms in the decompositions $n(b_0)$, $n(b_1)$, $n_T(f)$, $n_S(f)$ and $n(u_0)$ are fixed, while our parameters, namely: viscosity $\nu$, coefficients $b_{0m}$, $b_{1m}$, $f_m$ and $u_{0m}$, and $\left(A^{b_0}_l\right)_{l=1,\ldots,n(b_0)}$, $\left(A^{b_1}_l\right)_{l=1,\ldots,n(b_1)}$, $\left(A^f_{lp}\right)_{l=1,\ldots,n_T(f) ; p=1,\ldots,n_S(f)}$ and $\left(A^{u_0}_l\right)_{l=1,\ldots,n(u_0)}$ live in some Cartesian product of intervals $\mathcal P'$, subset of $\R^{1+4+n(b_0)+n(b_1)+n_T(f) n_S(f)+n(u_0)}$. Note that $m$ stands for \emph{mean} and is not a "numerical" index.

The compatibility condition \eqref{e:compatcond} constraints $b_{0m}$ and $b_{1m}$ as functions of the other parameters:
\[ b_{0m} = u_{0m} \;\;\text{ and }\;\; b_{1m}=u_{0m}+\sum_{l=1}^{n(u_0)} A_l^{u_0} \phiuz \]
so that our "compliant" parameters actually belong to $\mathcal P$ defined by:
\begin{multline}
\label{e:definitionP}
\mathcal P = \Big\{ \mu=\big(\nu,b_{0m},
							A^{b_0}_1,\ldots,A^{b_0}_{n(b_0)},
							b_{1m},A^{b_1}_1,\ldots,A^{b_0}_{n(b_1)}, 
							f_m, 
							A^f_{11},A^f_{12},\ldots,A^f_{1,n_S(f)}, \\
							A^f_{2,1},\ldots,A^f_{2,n_S(f)},\ldots,A^f_{n_T(f),n_S(f)}, 
							u_{0m},A^{u_0}_1,\ldots,A^{u_0}_{n(u_0)} \big) \in \mathcal P' 
							\text{ satisfying } \eqref{e:compatcond} \Big\}
\end{multline}

\subsection{Offline/online procedure}
\label{ss:redbcalc}
The key heuristic \cite{nguyen2005certified} for RB approximation of the (linear) parametrized variational problem is the following: given $\mu\in\mathcal P$, 
\[ \text{find } u(\mu)\in X \text{ so that } A(u(\mu),v;\mu)=L(v;\mu), \; \forall v \in X\]
is to choose a parameter-independent family $\mathcal R$ of linearly independent functions in $X$ --- with $\# \mathcal R \ll \dim X$, to achieve computational economy --- and then, given an instance of the parameter, to search for the reduced solution:
\[ \widetilde u(\mu) \in \text{Span}\mathcal R, \text{ so that } A(\widetilde u(\mu),v;\mu)=L(v;\mu) \; \forall v \in \text{Span}\mathcal R \]
Let us apply this idea on problem \eqref{e:disc}. We rely on the procedure described in \cite{nguyen2009reduced}, modified to allow parametrization of initial condition, boundary data and source term. To simplify notations, we do not explicitly write dependence of $u$ and $\widetilde u$ on $\mu$.

We suppose that a reduced basis $\mathcal R = \{\zeta_1, \ldots, \zeta_N \}$ has been chosen (see Subsection \ref{ss:choicebasis} for one way to do so); we define $\widetilde X = \text{Span}\mathcal R$ and we look for $\left\{\utk\right\}_{k=0,\ldots,\mathcal T} \subset \widetilde X$ satisfying:
\begin{subequations}
\label{e:discrb}
\begin{equation}
\label{e:discrb1}
\ut{0} = \widetilde \pi \left( \pi(u_0) \right) 
\end{equation}
and:
\begin{multline}
\label{e:discrb2}
\Big\langle \frac{\utk-\ut{k-1}}{\delt},v \Big\rangle  + c(\utk,\utk,v) + \nu a (\utk,v) + B(\utk,v) \\
\hspace*{2cm} =  \ell_\pi(v,t_k) + b_0(t_k) \beta_0(v) + b_1(t_k) \beta_1(v) \;\; \forall v \in \widetilde X \; \forall k=1,\ldots,\mathcal T
\end{multline}
\end{subequations}
where $\widetilde \pi$ is the orthogonal projection from $X$ onto $\widetilde X$, with respect to the standard $L^2$ inner product on $X$: $\langle w,v \rangle=\intzu wv$.

The offline/online procedure for computation of $\ut{0}$ will come easily from our parametrization of $u_0$ in Subsection \ref{ss:params}: since the constant function $\mathbf 1 = \sum_{j=0}^{\mathcal N} \phi_j$ belongs to $X$, we have:
\begin{equation}
\label{e:proju0}
\ut{0} = u_{0m} \widetilde\pi(\mathbf 1) + \sum_{l=1}^{n(u_0)} A^{u_0}_l \, \widetilde \pi \left( \pi \left( \phiuzp \right) \right)
\end{equation}

We now discuss computation of $\utk$ from $\ut{k-1}$ for $k=1,\ldots,\mathcal T$. We are willing to proceed with Newton steps as for the solution of \eqref{e:disc2}. We denote, for $k=1,\ldots,\mathcal T$, by
\[ \ut{k-1} = \sum_{j=1}^N u^{k-1}_j \zeta_j \quad;\quad \overline{\utk} = \sum_{j=1}^N \overline{u^k_j} \zeta_j \]
respectively, the reduced solution at time $t_{k-1}$, and previous guess for the reduced solution at time $t_k$.

Our procedure relies on the following proposition:
\begin{prpstn}
\begin{enumerate}
\item The Newton increment $\delta=\sum_{j=1}^N \delta_j \zeta_j$ satisfies the following equations:
\begin{multline}
\label{e:aaa}
\sum_{j=1}^N \delta_j \left\{ \frac{ \langle \zeta_j,\zeta_i \rangle }{\delt} + 2 \sum_{j'=1}^N \overline{u^k_{j'}} c(\zeta_{j'},\zeta_j,\zeta_i)
+ \nu a(\zeta_j,\zeta_i) + B(\zeta_j,\zeta_i) \right\} \\= 
\sum_{j=1}^N u_j^{k-1} \frac{ \langle \zeta_j, \zeta_i \rangle }{\delt} + \ell_\pi(\zeta_i,t_k) + b_0(t) \beta_0(\zeta_i) + b_1(t) \beta_1(\zeta_i)
\\- \sum_{j=1}^N \overline{u_j^k} \left( \frac{ \langle \zeta_j,\zeta_i \rangle }{\delt} + \sum_{j'=1}^N \overline{u_{j'}^k} c(\zeta_{j'},\zeta_j,\zeta_i)
+ \nu a (\zeta_j,\zeta_i) + B(\zeta_j,\zeta_i) \right) \;\; \forall i=1,\ldots,N
\end{multline}
\item We have:
\begin{equation}
\label{e:bbb}
\ell_\pi(\zeta_i,t_k) = f_m \int_0^1 \zeta_i + \sum_{l=1}^{n_T(f)} \sum_{p=1}^{n_S(f)} A_{lp}^f \phiftr(t_k) \intzu \pi\left( \phifsp \right) \zeta_i
\end{equation}
for all $i=1,\ldots,N$ and $k=1,\ldots,\mathcal T$.
\end{enumerate}
\end{prpstn}
\begin{proof} \begin{enumerate} 
\item Equation \eqref{e:discrb2} for $\utk=\overline{\utk}+\delta$ linearized near $\delta=0$, for $v=\zeta_i,\;\forall i=1,\ldots,N$ is:
\begin{multline*}
\frac{1}{\delt} \langle \overline{\utk}+\delta,\zeta_i \rangle + c(\overline{\utk},\overline{\utk},\zeta_i) + 2 c(\overline{\utk},\delta,\zeta_i) + \nu a (\overline{\utk}+\delta,\zeta_i) + B(\overline{\utk}+\delta,\zeta_i) \\ = \frac{1}{\delt} \langle \ut{k-1},\zeta_i \rangle + \ell_\pi(\zeta_i,t_k) + b_0(t) \beta_0(\zeta_i) + b_1(t) \beta_1(\zeta_i) \;\; \forall i=1,\ldots, N 
\end{multline*}
Rewriting this equation using expansions of $\overline{\utk}$ and $\delta$ in $\mathcal R$ and linearity of $\langle \cdot,\cdot \rangle$, $c$, $a$ and $B$ with respect to their first argument, and putting all $(\delta_j)_j$-dependent terms on the left-hand side give the announced equation.
\item is a direct consequence of the parametrization of $f$ given in Subsection \ref{ss:params}.  
\end{enumerate}\end{proof}

The following Proposition \ref{prop:welldefinednewton} justifies the well-definedness of the Newton iteration, for an orthonormal reduced basis $\{\zeta_1,\ldots,\zeta_N\}$. In practice, the Gram-Schmidt process can always be used to ensure this condition.

\begin{prpstn}
\label{prop:welldefinednewton}
	Suppose that $\{\zeta_1,\ldots,\zeta_N\}$ is orthonormal with respect to $\langle \cdot , \cdot \rangle$.

	Then, for $\Delta t$ small enough, i.e. $\Delta t < \Delta t^*\left(\nu,\{\zeta_i\}_{i=1,\ldots,N}\right)$, and initial guess $\overline{\widetilde u^k}$ sufficiently close to $\widetilde u^k$, the Newton iteration \eqref{e:aaa} is well defined and converges (quadratically) to $\widetilde u^k$.
\end{prpstn}
\begin{proof}
	Iteration \eqref{e:aaa} is a Newton iteration for solving $F(x)=\alpha$, for appropriate $\alpha$ and $F$ given by:
	\[ F(x)=\sum_{i=1}^N \left( \frac{\langle x,\zeta_i \rangle}{\Delta t}+c(x,x,\zeta_i)+\nu a(x,\zeta_i)+B(x,\zeta_i) \right) \zeta_i \]
	We apply the result stated pp. 353--355 of \cite{quarteroni2008numerical} and pp. 362--367 of \cite{girault1986finite}. To do so, we have to check that, for $\Delta t$ small enough:
	\begin{enumerate}
		\item the differential of $F$ at $\widetilde u^k$, denoted by $DF(\widetilde u^k)$, is invertible; 
		\item $DF(\cdot)$ is Lipschitz-continuous \emph{i.e.} there exist $L>0$ so that 
			\[ \forall x,x' \in X,\;\;\; \forall v \in X, \;\;\; \norm{ DF(x)\cdot v-DF(x')\cdot v } \leq L \norm{x-x'} \norm{v} \]
	\end{enumerate}
	It is easy to check that the matrix of $DF(\widetilde u^k)$ in the reduced basis $\{\zeta_1,\ldots,\zeta_N\}$ is diagonally dominant, hence invertible, for:
	\[ \Delta t <  \Delta t^* := \min_{i=1,\ldots,N} \frac{1}{\sum_{j=1}^N \abs{2c(\widetilde u^k,\zeta_j,\zeta_i)+\nu a(\zeta_j,\zeta_i)+B(\zeta_j,\zeta_i)}} \]
	and (2) is a straightforward computation.
\end{proof}

We put our offline/online procedure in Algorithm \ref{algo:1}.
\begin{algo}
\caption{offline/online procedure}
\label{algo:1}
\begin{itemize}
\item \emph{offline:} 
\begin{enumerate}
 \item choose a parameter-, and time-independent reduced basis $\{\zeta_1,\ldots,\zeta_N\}$ (see Section \ref{ss:choicebasis})
 \item compute and store the following parameter-independent functions of $\widetilde X$ and scalars, for all $i,j,j'=1,\ldots,N$, $l=1,\ldots,n(u_0)$, $p=1,\ldots,n_S(f)$:
\[ \begin{array}{cc}
\widetilde\pi(\mathbf 1) &
\widetilde \pi \left( \pi \left( \phiuzp \right) \right) \\
\langle \zeta_j,\zeta_i \rangle &
a(\zeta_j,\zeta_i) \\
c(\zeta_{j'},\zeta_j,\zeta_i) &
B(\zeta_j,\zeta_i) \\
\beta_0(\zeta_i) &
\beta_1(\zeta_i) \\
\int_0^1 \zeta_i &
\intzu \pi\left( \phifsp \right) \zeta_i \\
\end{array} \]
\end{enumerate}
\item \emph{online:} 
  \begin{enumerate}
  \item assemble $\ut{0}$ as the linear combination \eqref{e:proju0} ;
  \item for $k=1, \ldots, \mathcal T$:
  \begin{enumerate}
  		\item set up an initial guess $\overline{\utk} = \sum_{j=1}^N \overline{u^k_j} \zeta_j$ ;
		\item compute and store $\ell_\pi(\zeta_i,t_k)$ by using \eqref{e:bbb} ;
		\item look for $\delta=\sum_{j=1}^N \delta_j \zeta_j$ by solving the linear system \eqref{e:aaa} ;
		\item set $\utk \leftarrow \overline{\utk}+\delta$ ;
		\item if $\norm{\delta}$ is small enough:
		  \begin{enumerate}
		   \item output $\utk$ 
			\item set $k=k+1$
			\end{enumerate}
		 \item else:
		  \begin{enumerate}
		   \item update the guess : \mbox{$\overline{\utk} \leftarrow \utk$, i.e. ${\utk}_j \leftarrow u^k_j \; \forall j=1,\ldots,N$}
		   \item go back to (c)
		  \end{enumerate}
	\end{enumerate}
	\end{enumerate}
\end{itemize}
\end{algo}
Let us make some remarks about the complexity of the above online algorithm, in contrast with the "full basis" one described in Section \ref{ss:num}:
\begin{rmrk}
	\label{rmrk:a}
\begin{enumerate}
\item The most computationally demanding step is (2) (c), since it involves resolution of a (nonsymmetric, dense) $N \times N$ linear system; the "full basis" counterpart solves $(\mathcal N + 1) \times (\mathcal N+1)$ nonsymmetric tridiagonal system. Thus significant computational savings are expected for $N \ll \mathcal N$.
\item Thanks to our parametrization of $f(t,\cdot)$, all integrals over $[0;1]$ in equation \eqref{e:aaa} can be precomputed during the offline phase, yielding a $\mathcal N$-\emph{independent} online phase. This means that one can in principle choose \emph{arbitrary} high precision on the full model \emph{without} impact on the marginal cost of evaluation of an online solution. This "$\mathcal N$-independence" property shall be required of every complexity of any online procedure.
One should also note that our online procedure does not produce "nodal" values $\utk(x_i), \, i=0,\ldots,\mathcal N$ (as this would clearly violate the $\mathcal N$-independence), but rather the components of $\utk$ in the reduced basis.
\item Taken independently, the number of parameters $n(u_0)$, $n(b_0)$, $n(b_1)$, $n_S(f)$ and $n_T(f)$, as well as the number of timesteps $\mathcal T$ have a linear impact on the online complexity. Moreover, due to the double sum in \eqref{e:bbb}, the online complexity is proportional to $n_S(f) n_T(f)$. We note that an advantage of treating the Dirichlet boundary condition weakly, as we do in this paper, is that $n_T(f)$ does not get increased by functions of $n({b_0})$ or $n({b_1})$. This is a clear advantage of our method: when boundary conditions are treated by returning to an homogeneous Dirichlet problem, as in \cite{jung2009reduced}, $n(b_0)+n(b_1)$ terms are added in the parametrization of $f$.

\end{enumerate}
\end{rmrk}

\subsection{Choice of the reduced basis}
\label{ss:choicebasis}
In this subsection, we describe different ways of choosing a pertinent reduced basis $\{\zeta_1,\ldots,\zeta_N\}$. These lead to three different bases fitting into the certified (that is to say, the three admit the same procedure for online error bound we describe in Section \ref{s:error}) reduced basis framework.

The first is based on proper orthogonal decomposition (POD), the second is based on a "greedy" selection algorithm. The third is an hybridation of POD and greedy. These methods are standard in the literature.

\subsubsection*{Notation}
The two procedures described below involve computation of the reference solution for different instances of the parameter, so we should use special notations, local to this section, to emphasize the dependence of the reference solution on the parameters. We define a parametrized solution $u$ by:
\[
u : \left\{ \begin{aligned} \{1,\ldots,\mathcal T\} \times \mathcal P &\rightarrow X \\
\left( k, \mu \right) &\mapsto u(k,\mu) = \utk \text{ satisfying \eqref{e:disc2} for } \mu \text{ as parameter} \end{aligned} \right.
\]

\subsubsection{POD-driven procedure}
We denote by $\{u_j^k(\mu)\}_j$ the coordinates of $u(k,\mu)$ in the basis $\left\{\phi_0,\ldots,\phi_{\mathcal N}\right\}$:
\[ u(k,\mu) = \sum_{j=1}^{\mathcal N} u_j^k(\mu) \phi_j \]

In the POD-based procedure (see \cite{chatterjee2000introduction}) of the reduced basis choice, we choose a finite-sized parameter sample $\Xi \subset \mathcal P$, compute the reference solutions $u(k,\mu)$ for all $k=1,\ldots,\mathcal T$ and all $\mu\in\Xi$, and form the \emph{snapshots matrix} containing the components of these solutions in our basis $\left\{\phi_0,\ldots,\phi_{\mathcal N}\right\}$ :
\[
M=\begin{pmatrix}
	u_0^0(\mu_1) & u_0^1(\mu_1) & \cdots & u_0^{\mathcal T}(\mu_1) & u_0^0(\mu_2) & \cdots & \cdots & u_0^{\mathcal T}(\mu_S) \\
	u_1^0(\mu_1) & u_1^1(\mu_1) & \cdots & u_1^{\mathcal T}(\mu_1) & u_1^0(\mu_2) & \cdots & \cdots & u_1^{\mathcal T}(\mu_S) \\
	\vdots		 & \vdots		 & \vdots & \vdots						& \vdots			& \vdots & \vdots & \vdots						  \\
	u_{\mathcal N}^0(\mu_1) & u_{\mathcal N}^1(\mu_1) & \cdots & u_{\mathcal N}^{\mathcal T}(\mu_1) & u_{\mathcal N}^0(\mu_2) & \cdots & \cdots & u_{\mathcal N}^{\mathcal T}(\mu_S) 
	\end{pmatrix} \]
where $\Xi = \{ \mu_1, \ldots, \mu_S \}$.

One can check that $M$ has $\mathcal N+1$ rows and $S(\mathcal T+1)$ columns.

To finish, we choose the size $N < S(\mathcal T+1)$ of the desired reduced basis, we form the $S(\mathcal T+1) \times S(\mathcal T+1)$ non-negative symmetric matrix $M^T \Omega M$, where $\Omega$ is the matrix of our inner product $\langle \cdot , \cdot \rangle$, to find $z_1, \ldots, z_N$ the $N$ leading nonzero eigenvectors of this matrix (that is, the ones associated with the $N$ largest eigenvalues, counting repeatedly possible nonsimple eigenvalues), and, for $i=1,\ldots,N$, the coordinates of $\zeta_i$ with respect to $\{\phi_0,\ldots,\phi_{\mathcal N}\}$ are given by:
\begin{equation}
\label{e:podeqn}
 \frac{1}{\norm {Mz_i}} M z_i  
\end{equation}

\subsubsection{"Local" greedy selection procedure}
\label{sss:greedy}
The \emph{greedy} basis selection algorithm (cf. \cite{grepl2005posteriori,grepl2005thesis}) is the following: 
\begin{algo}
\caption{Greedy basis selection}
Parameter: $N$, the desired size of the reduced basis.
\begin{enumerate}
	\item Choose a finite-sized, random, large sample of parameters $\Xi\subset\mathcal P$.
	\item Choose $\mu_1\in\mathcal{P}$ and $k_1\in\{0,\ldots,\mathcal T\}$ at random, and set 
	\[ \zeta_1=\frac{u(\mu_1, k_1)}{\norm{u(\mu_1,k_1)}} \]
	\item Repeat, for $n$ from $2$ to $N$:
	\begin{enumerate}
		\item Find:
			\[ (\mu_n,k_n) = \argmax_{(\mu,k)\in\Xi\times\{0,\ldots,\mathcal T\}} \epsilon^*(\mu,k) \]
			where $\epsilon^*(\mu,k)$ is a (fastly evaluated) estimator of the RB error $\norm{u(\mu,k) - \widetilde u(\mu,k)}$, where $\widetilde u(\mu,k)$ stands for the RB approximation to $u(\mu,k)$ using $\{\zeta_1,\ldots,\zeta_{n-1}\}$ as reduced basis (see below).
		\item Compute $\zeta_n^* = u(\mu_n,k_n)$.
		\item Using one step of the (stabilized) Gram-Schmidt algorithm, find $\zeta_n \in \textrm{Span}\{\zeta_1,\ldots,\zeta_{n-1},\zeta_n^*\}$ so that $\{\zeta_1,\ldots,\zeta_{n-1},\zeta_n\}$ is an orthonormal family of $(X,\langle \cdot, \cdot \rangle)$.
	\end{enumerate}
\end{enumerate}
\end{algo}

The "greedy" name for this algorithm comes from the fact that the algorithm chooses, at each step of the repeat loop, the "best possible" time and parameter tuple to the reduced basis, that is the one for which the RB approximation error is estimated to be the worst.

Let's now discuss the choice for the error indicator $\epsilon^*$. A natural candidate would be the online error bound $\epsilon$ described in Section \ref{s:error}. However, as we shall see in the next section, the bound for timestep $t_k$ is a compound of the "propagation" of the error made in the previous timesteps $t_0, t_1, \ldots, t_{k-1}$ (the $\epsilon_{k-1}$ term) and the "local error" just introduced at the $k$-th time discretization. Hence, this error estimator $\epsilon_k$ has a natural tendency to grow (exponentially) with $k$. Thus using it as error indicator $\epsilon^*$  will favor times $k_n$ near final time $\mathcal T$ to be chosen at step 3.(a) of the algorithm below. Such choices are suboptimal, because including them in the reduced basis will not fix this exponential growth problem which is inherent to our approximation. Instead we use, as in \cite{grepl2005posteriori}, a purely \emph{local}-at-$t_k$ error indicator, that is: the computable error bound described in Section \ref{s:error} when taking $\epsilon_{k-1}=0$. 

It has been noted in the literature that the greedy procedure can ``stall'' (i.e. select vectors for addition in the basis which have no effect in decreasing the error bounds). The author in \cite{grepl2005thesis} provides a ``back-up procedure'' in such a case. In our experiments, however, we have not noticed any ``stalling'' of the greedy procedure, maybe because of our sharper error bound.

\subsubsection{POD-Greedy procedure}
We can also make use of the POD-Greedy procedure \cite{haasdonk2008reduced} in Algorithm \ref{algo:podgdy}.
This procedure aims to combine the advantages of the greedy and the POD based procedure; it has also been proposed so as to overcome the ``stalling'' of the greedy procedure. 
\begin{algo}
	\caption{POD-Greedy basis selection}
	\label{algo:podgdy}
Parameter: $P_1$.
\begin{enumerate}
	\item Choose a finite-sized, random, large sample of parameters $\Xi\subset\mathcal P$.
	\item Choose $\mu_1\in\mathcal{P}$ at random, and choose as the current reduced basis $\mathcal B$ an orthonormalized basis of $\text{Span} \{ u(\mu_1,t_0), u(\mu_1,t_1), \ldots, u(\mu_1, \mathcal T) \}$.
	\item Repeat, until the desired number of items in the basis is reached:
	\begin{enumerate}
		\item Find:
			\[ \mu^* = \argmax_{\mu\in\Xi} \epsilon^*(\mu) \]
			where $\epsilon^*(\mu)$ is a (fastly evaluated) estimator of the RB error $\norm{u(\mu) - \widetilde u(\mu)}$, where $\widetilde u(\mu)$ stands for the RB approximation to $u(\mu)$ using $\mathcal B$ as reduced basis.
		\item Append to $\mathcal B$ the $P_1$ leading POD modes of $\{ u^{\text{proj}}(\mu^*,t_0), u^{\text{proj}}(\mu^*,t_1), \ldots, u^{\text{proj}}(\mu^*, \mathcal T) \}$, where the proj superscript denotes projection on the orthogonal complement of $\text{Span}(\mathcal B)$.
	\end{enumerate}
        \item Output $\mathcal B=\{\zeta_1,\ldots,\zeta_{\#\mathcal B}\}$ as reduced basis.
\end{enumerate}
\end{algo}
Again, an error indicator $\epsilon^*$ has to be chosen. We do as in \cite{knezevic2010certified}, where the authors take the online error bound at final time.

\subsubsection{Expansion of the basis by initial data modes}
\label{sss:enrich}
In case they did not get chosen by the POD or greedy algorithm, a classical strategy is to initialize the basis selection algorithms by taking in the reduced basis the constant function $\mathbf 1$, and the functions $\phiuzp$ for $l=1,\ldots,n(u_0)$. This may increase the size of the reduced basis (and thus online computation times) but ensures that $\ut{0} = u_0$ (i.e. initial error is zero). 
Such an enrichment can possibly be a good move, as the error gets accumulated and amplified throughout the time iterations, zero initial error will certainly reduce the (estimated, as well as actual) RB approximation error.

\section{Error bound}
\label{s:error}
In this section, we derive a parameter and time dependent online error bound $\epsilon^k$ (for $k=1,\ldots,\mathcal T$) satisfying:
$ \norm{u^k_e - \utk} \leq \epsilon^k, $
where $\norm{v}=\left( \intzu v^2 \right)^{1/2}$, and $u^k_e$ is our reference ``truth'' solution in $X$ satisfying the fully discretized PDE with strong Dirichlet enforcement, that is:
\begin{equation}
\label{e:uktruth}
\left\{ \begin{array}{l}
\frac{ \langle u_e^k - u_e^{k-1}, v \rangle }{\Delta t} + c (u^k_e,u^k_e,v) + \nu a (u^k_e,v) = \ell_\pi(v,t_k)\;\;\; \forall v \in X_0 \\
u^k_e(0)=b_0(t_k) \\
u^k_e(1)=b_1(t_k) \end{array} \right.
\end{equation}
where $X_0$ is the ``homogeneous'' subspace of $X$:
\[ X_0 = \{ v \in X \text{ st. } v(0)=v(1)=0 \}. \]

Our error bound should be precise enough (i.e., not overestimating the actual error $\norm{u^k_e - \utk}$ too much, and approaching zero as $N$ increases) and online-efficient (that is, admit an offline/online computation procedure with an $\mathcal N$-independent online complexity).

We notice that our error bound measures the error between the reduced and the reference solution; it does not reflect the discretization error made when replacing the actual analytical solution of the Burgers equation with its numerical approximation \eqref{e:uktruth}. This is consistent with the fact that RB methods relies strongly on the existence of a high-fidelity numerical approximation of the analytical solution by a discrete one, hence regarded as ``truth''.

Subsection \ref{ss:ebound} deals with the derivation of the error bound; this error involves quantities whose computation is detailed in Subsections \ref{ss:initerr} and \ref{ss:znres}.

\subsection{Error bound}
\label{ss:ebound}
The sketch of the derivation of our error bound is the following: we first give a "theoretical" error bound in Theorem \ref{thm:1}; we then replace the uncomputable quantities appearing in this bound by their computable surrogates in paragraph \ref{sss:compt}.

\subsubsection{Theoretical error bound}
\textbf{Notation. }
We suppose that the convergence tests appearing in the Newton iterations performed in Section \ref{ss:numT} and Section \ref{ss:redbcalc} are sufficiently demanding so as to neglect the errors due to the iterative solution of the nonlinear systems \eqref{e:nleq} and \eqref{e:discrb2}.

We now set up some notations : first the error at time $t_k$:
\[ e_k = \begin{cases} u^k_e - \utk &\text{ if } k>0 \\
							\pi(u_0) - \widetilde \pi \left( \pi(u_0) \right) &\text{ if } k=0\end{cases}	\]
the residual form $r_k$, for $v \in X_0$:
\begin{equation}
\label{e:resform}
r_k(v) = \ell_\pi(v,t_k)-\frac{1}{\Delta t} \langle \utk - \ut{k-1}, v \rangle
- c(\utk, \utk, v)-\nu a(\utk,v) 
\end{equation}
and the "$X_0$-norm" of the residual:
\begin{equation*}
	\norm{r_k}_0 = \sup_{v\in X_0, \norm{v}=1} r_k(v)
\end{equation*}

We introduce :
\[ \psi_k(v,w) = 2c(\utk,v,w)+\nu a(v,w) \]
and the so-called \emph{stability constants}:
\begin{equation}
\label{e:defCk}
C_k= \inf_{v\in X_0, \norm{v}=1} \psi_k(v,v) 
\end{equation}
To finish, we define:
\[ \eta_k = \abs{e_k(0)} \norm{\phi_0} + \abs{e_k(1)} \norm{\phi_{\mathcal N}} \quad;\quad 
 \sigma_k = 2 \abs{C_k} \eta_k \quad;\quad
 \mathcal E = \sup_{v\in X_0, \norm v=1} v\left( \frac{1}{\mathcal N} \right) \]
($\mathcal E$ is finite because $X_0$ is finite dimensional),
and, finally:
\[ f_k = \mathcal E \big( \abs{e_k(0)} \abs{ \psi_k(\phi_0,\phi_1) } +
\abs{e_k(1)} \abs{ \psi_k(\phi_{\mathcal N},\phi_{\mathcal N-1}) } \big) \]
\[ \xi_k^{\mathcal A} = \frac{ \mathcal E^2(\abs{e_k(0)}+\abs{e_k(1)})}{3} \quad;\quad 
	\xi_k^{\mathcal B} = \frac{5}{3} \mathcal E \left( e_k(0)^2+e_k(1)^2 \right) \quad;\quad 
   \xi_k^{\gamma} = \frac{ \abs{e_k(0)}^3+\abs{e_k(1)}^3 }{3}  
\]
We will also make use of the standard notation:
\[ \partieneg{x} = \max(-x, 0) \;\;\; \forall x \in \R \]
The theoretical foundation for our error bound is the following theorem. As the technique developed in \cite{nguyen2009reduced} did not fit our problem (because of the weak Dirichlet treatment, whose advantage has been shown in Remark \ref{rmrk:a} point (3)), we developed a new strategy for obtaining this error bound. 
\begin{thrm}
\label{thm:1}
If:
\begin{equation}
\label{e:hypCn}
\frac{1}{\delt} + C_k - \xi_k^{\mathcal A} > 0 \; \forall k=1,\ldots,\mathcal T
\end{equation}
then the norm of the error $\norm{e_k}$ satisfies:
\[ \norm{e_k} \leq \begin{cases}
		\frac{\mathcal B_k + \sqrt {\mathcal D_k}}{2 \mathcal A_k} & \text{ if } \mathcal D_k\geq 0 \\
		\frac{\mathcal B_k}{\mathcal A_k} & \text{ if } \mathcal D_k< 0 \end{cases} \]
with:
\[ \mathcal{A}_k = \frac{1}{\delt} + C_k - \xi_k^{\mathcal A} \quad;\quad
 \mathcal{B}_k = \frac{2 \eta_k + \norm{e_{k-1}} + \mathcal E \langle \phi_0,\phi_1 \rangle ( \abs{e_k(0)}+\abs{e_k(1)} )}{\delt}  + \sigma_k + f_k + \norm{r_k}_0 + \xi_k^{\mathcal B} \]
\begin{multline*}
	\gamma_k=  \frac{\eta_k \norm{e_{k-1}} +\mathcal E \eta_k \langle \phi_0,\phi_1 \rangle (\abs{e_k(0)}+\abs{e_k(1)})  }{\Delta t} + \eta_k f_k + \partieneg{C_k} \eta_k^2 + \frac{1}{6}\abs{e_k(1)^3-e_k(0)^3} + \xi_k^\gamma + \norm{r_k}_0\eta_k
\end{multline*}
and:
\[ \mathcal D_k = \left(\mathcal B_k\right)^2 + 4 \mathcal A_k \gamma_k. \]
\end{thrm}

The proof of Theorem \ref{thm:1} is presented in the appendix.

\subsubsection{Computable error bound. }
\label{sss:compt}
We now find an efficiently computable (that is, with an offline/online decomposition, with a complexity of the online part independent of $\mathcal N$) error bound $\epsilon_k$ derived from the one described above; to do so we discuss each of the ingredients appearing in its expression. 
\begin{itemize}
\item Computation of the norm of the initial error $\norm{e_0}$ is addressed in the next Section \ref{ss:initerr}; the one of $\norm{r_k}_0$ is in Section \ref{ss:znresidual}. 
\item We have, for $w\in\{0,1\}$:
\[ e_k(w) = u^k_e(w)-\utk(w) = b_w(t_k)-\utk(w) \]
so that $e_k(0)$ and $e_k(1)$ can be computed during the online phase.
\item 
	Similarly, the scalars $\norm{\phi_0}$, $\psi_k(\phi_0,\phi_1)$ and $\psi_k(\phi_{\mathcal N,\mathcal N-1})$  can straightforwardly be computed online.
\item The "continuity constant" $\mathcal E$ can be computed offline and stored by solving the optimization problem defining it.
 Thus $\eta_k$, $f_k$, $\xi_k^{\mathcal A}$, $\xi_k^{\mathcal B}$ and $\xi_k^\gamma$ can be computed online.
\item The exact value of $C_k$ could be found by solving a generalized eigenvalue problem on $X$: $C_k$ is the smallest $\lambda\in\R$ so that there exists $z\in X_0$, $\norm z=1$ satisfying:
\[ \psi_k^{Sym}(z,v) = \lambda \langle z,v \rangle \;\; \forall v \in X_0 \]
with $\psi_k^{Sym}$ the symmetric bilinear form defined by:
\[ \psi_k^{Sym}(w,v) = \nu a(w,v) + c(\utk,w,v) + c(\utk,v,w) \;\; \forall w,v\in X_0 \]
The cost of doing so is prohibitive as it is an increasing function of $\dim X = \mathcal N+1$. Instead, we will see in Section \ref{ss:scm} how to compute lower and upper bounds $C_k^{inf}$ and $C_k^{sup}$:
$ C_k^{inf} \leq C_k \leq C_k^{sup} $.
\item We can then compute the following lower and upper bounds for $\mathcal A_k$:
	\[ \mathcal A_k^{inf} = \frac{1}{\delt} + C_k^{inf} - \xi_k^{\mathcal A} \quad\quad;\quad\quad 
		\mathcal A_k^{sup} = \frac{1}{\delt} + C_k^{sup} - \xi_k^{\mathcal A} \]
and the hypothesis \eqref{e:hypCn} is ensured by checking that $\mathcal A_k^{inf} > 0$.
\item We can also compute an upper bound of $\sigma_k$:
\[ \sigma_k^{sup} = 2 \eta_k \max\big( \abs{C_k^{sup}}, \abs{C_k^{inf}} \big) \]
\item To compute an upper bound of $\mathcal B_k$, we need to replace the preceding error norm $\norm{e_{k-1}}$ which is (except for $k=1$) not exactly computable, with the online upper bound $\epsilon_{k-1} \geq \norm{e_{k-1}}$ at the preceding time step:
	\[ \mathcal B_k^{sup} = \frac{2 \eta_k + \epsilon_{k-1} + \mathcal E \langle \phi_0,\phi_1 \rangle ( \abs{e_k(0)}+\abs{e_k(1)} )}{\delt}  + \sigma_k^{sup} + f_k + \norm{r_k}_0 + \xi_k^{\mathcal B}  \]
\item And $\gamma_k$ gets replaced by its upper bound $\gamma_k^{sup}$:
\begin{multline*}
	\gamma_k^{sup}= \frac{\eta_k \epsilon_{k-1} +\mathcal E \eta_k \langle \phi_0,\phi_1 \rangle (\abs{e_k(0)}+\abs{e_k(1)})  }{\Delta t} + \eta_k f_k + \partieneg{C_k^{inf}} \eta_k^2  + \frac{1}{6}\abs{e_k(1)^3-e_k(0)^3} + \xi_k^\gamma + \norm{r_k}_0\eta_k
\end{multline*}
\item We finally compute an upper bound for $\mathcal D_k$:
	\[ \mathcal D_k^{sup} = \begin{cases} \left(\mathcal B_k^{sup} \right)^2 + 4 \mathcal A_k^{sup} \gamma_k^{sup} &\text{ if } \gamma_k^{sup}\geq 0 \\
\left(\mathcal B_k^{sup} \right)^2 + 4 \mathcal A_k^{inf} \gamma_k^{sup} &\text{ if } \gamma_k^{sup}< 0. \end{cases} \]
and our "computable" error bound is then:
\[ \begin{cases}
		\frac{\mathcal B_k^{sup} + \sqrt{ \mathcal D_k^{sup} } }{2 \mathcal A_k^{inf}} & \text{ if } \mathcal D_k^{sup} \geq 0 \\
		\frac{\mathcal B_k^{sup}}{\mathcal A_k^{inf}} & \text{ if } \mathcal D_k^{sup} < 0.
	\end{cases}
	\]
\end{itemize}

The remainder of the section consists in the description of computation of the four left-out quantities $\norm{e_0}$, $\norm{r_k}_0$, $C_k^{sup}$ and $C_k^{inf}$.

\subsection{Initial error}
The present subsection deals with efficient computation of the $\norm{e_0}$ term in the computable error bound described in paragraph \ref{sss:compt}.
\label{ss:initerr}
We denote by  $\mathbf{H}$ the Gram matrix of the family:
\[ \Big\{ \mathbf{1}-\widetilde\pi(\mathbf 1), \pi\left(\Phi^{u_0}_1(\cdot) \right) - \widetilde\pi\left(\pi\left( \Phi^{u_0}_1(\cdot) \right)\right), \ldots, 
		\pi\left( \Phi^{u_0}_{n(u_0)}(\cdot) \right) - \widetilde\pi\left(\pi\left(\Phi^{u_0}_{n(u_0)}(\cdot) \right)\right) \Big\} 
\]
that is, $\mathbf{H}$ is the $(1+n(u_0))\times(1+n(u_0))$ symmetric matrix of all the inner products between two any of the above vectors 
and by $\mathbf{e_0}$ the vector containing the components of $e_0$ with respect to the family above, \emph{i.e.}:
\[ \mathbf{e_0} = \left( u_{0m}, A_1^{u_0}, \ldots, A_{n(u_0)}^{u_0} \right)^T \]

We have:
\begin{lmm}
The norm of the initial error $\norm{e_0}$ is given by:
\begin{equation}\label{e:e0} ||e_0|| = \sqrt{ \mathbf{e_0}^T \mathbf{H} \mathbf{e_0} } \end{equation}
\end{lmm}
\begin{proof} Parametrization \ref{ss:params} gives:
\begin{align*}
e_0 &= \pi(u_0) - \widetilde\pi\left( \pi(u_0) \right) 
= u_{0m} \left( \mathbf{1}-\widetilde\pi(\mathbf 1) \right) + \sum_{l=1}^{n(u_0)} A_l^{u_0} 
\left( \pi\left(\Phi^{u_0}_l( \cdot)\right) - \widetilde\pi\left(\pi\left(\Phi^{u_0}_l( \cdot )\right) \right) \right)
\end{align*}
and the result follows from the expansion of $\norm{e_0}^2$ when $e_0$ is replaced by the expression above. \end{proof}

This formula allows us to compute the time and parameter-independent Gram matrix $\mathbf{H}$ during the offline phase, and, during the online phase, to assemble the $(1+n(u_0))$-vector $\mathbf{e_0}$ and to perform \eqref{e:e0} to get $||e_0||$ with an online cost dependent only of $n(u_0)$.

\subsection{Norm of the residual}
We now present the computation of the $\norm{r_k}_0$ term in the computable error bound of Section \ref{sss:compt}.

\label{ss:znres}
\begin{itemize}
\item Let the expansions of the reduced solutions with respect to the reduced basis be:
\[ \ut{p} = \sum_{j=1}^N u^p_j \zeta_j \;\;\; \text{ for } p\in\{k-1, k\} \] 
\item Let $\mathbf G$ be the $(1+n_S(f)+2N+N^2) \times (1+n_S(f)+2N+N^2)$-sized Gram matrix of 
\begin{multline*}
\{\Gamma^{int}, \Gamma^{fS}_1, \ldots, \Gamma^{fS}_{n_S(f)}, 
\Gamma^{\langle\rangle}_1, \ldots, \Gamma^{\langle  \rangle}_N, 
\Gamma^c_{1,1}, \Gamma^c_{1,2}, \ldots, \Gamma^c_{1,N}, 
\Gamma^c_{2,1}, \ldots, \Gamma^c_{2,N}, \ldots, \Gamma^c_{N,N},
\Gamma^a_1, \ldots, \Gamma^a_N \}
\end{multline*}
where $\Gamma^{int}, \Gamma^{fS}_p, \Gamma^{\langle  \rangle}_j, \Gamma^c_{j,j'}, \Gamma^a_j \in X_0$ ($p=1,\ldots,n_S(f); j,j'=1,\ldots,N$) satisfy:
\begin{align*}
	\langle \Gamma^{int}, v \rangle &= \intzu v &\;\;\; \forall v \in X_0 \\
	\langle \Gamma^{fS}_p, v \rangle &= \intzu \pi\left( \phifsp \right) v  &\;\;\; \forall v \in X_0 \\
	\langle \Gamma^{\langle\rangle}_j, v \rangle &= \langle\zeta_j,v\rangle &\;\;\; \forall v \in X_0 \\
	\langle \Gamma^{c}_{j,j'}, v \rangle &= c(\zeta_j,\zeta_{j'},v) &\;\;\; \forall v \in X_0 \\
	\langle \Gamma^{a}_j, v \rangle &= a(\zeta_j,v) &\;\;\; \forall v \in X_0 
\end{align*}
(Those $\Gamma$'s exist by virtue of the Riesz representation theorem).
\item Let $\rhog_k$ be the following vector:
\begin{multline*}
\rhog_k = \Big(
f_m, 
\sum_{l=1}^{n_T(f)} A_{l,j}^f \phiftr(t_k)\; (j=1,\ldots,n_S(f)), 
- \frac{1}{\delt} (u^k_j - u^{k-1}_j)\; (j=1,\ldots,N), \\
u^k_j u^k_l\; (j,l=1,\ldots,N) ,
\nu u^k_j\; (j=1,\ldots,N) 
\Big)^T
\end{multline*}
\end{itemize}

\label{ss:znresidual}
\begin{lmm}
We have:
\begin{equation}
\label{e:residu0}
\norm{r_k}_0 = \sqrt{ \rhog_k^T \mathbf{G} \rhog_k } 
\end{equation}
\end{lmm}

\begin{proof}
From the Riesz representation theorem, there exists a unique $\rho_k \in X_0$ so that
\begin{equation}
\label{e:dualrk}
\langle \rho_k, v \rangle = r_k(v) \;\; \forall v \in X_0 
\end{equation}
and we have:
$\norm{r_k}_0 = \norm{\rho_k}$.

From \eqref{e:dualrk} and the definition of $r_k$ \eqref{e:resform}, we have that $\rho_k$ is defined uniquely by:
\[ \langle \rho_k,v \rangle = \ell_\pi(v, t_k)-\frac{1}{\Delta t} \langle \utk - \ut{k-1}, v \rangle
- c(\utk, \utk, v)-\nu a(\utk,v) \forall v \in X_0 \]
because $\beta_0(v)=\beta_1(v)=B(\cdot,v)=0$ for all $v\in X_0$.

Using parametrization \eqref{e:bbb} of  $\ell_\pi(\cdot,t_k)$, we get that \eqref{e:dualrk} is equivalent to:
\begin{multline}
\label{e:dualrk2}
\langle \rho_k,v \rangle = f_m \int_0^1 v + \sum_{l=1}^{n_T(f)} \sum_{p=1}^{n_S(f)} A_{lp}^f \phiftrt(t_k) \intzu \pi\left( \phifsp \right) v 
 - \frac{1}{\delt} \sum_{j=1}^N (u^k_j - u^{k-1}_j) \langle \zeta_j,v \rangle 
\\ - \sum_{j=1}^N u^k_j \left( \sum_{j'=1}^N u^k_{j'} c(\zeta_j,\zeta_{j'},v) + \nu a(\zeta_j,v) \right) \;\forall v\in X_0 
\end{multline}
By the superposition principle, $\rho_k$ can be written as the linear combination:
\begin{multline*}
\rho_k = f_m \Gamma^{int} + \sum_{p=1}^{n_S(f)} \left( \sum_{l=1}^{n_T(f)} A_{lp}^f \phiftr(t_k) \right) \Gamma^{fS}_p
  - \frac{1}{\delt} \sum_{j=1}^N (u^k_j - u^{k-1}_j) \Gamma^{\langle\rangle}_j
 - \sum_{j=1}^N \sum_{j'=1}^N u^k_j u^k_{j'} \Gamma^c_{j,j'}
 - \sum_{j=1}^N \nu u^k_j \Gamma^a_j
\end{multline*}
Thus $\rhog_k$ contains the components of $\rho_k$ with respect to the family whose $\mathbf G$ is the Gram matrix, and so:
\begin{equation}
\norm{r_k}_0 = \norm{\rho_k} = \sqrt{ \rhog_k^T \mathbf{G} \rhog_k }.  
\end{equation}
\end{proof}

The offline/online decomposition for computation of $\norm{r_k}_0$ is as follows: in the offline phase, we compute the $\Gamma$'s vectors, and compute and store their Gram matrix $\mathbf G$. In the online phase, we compute $\rhog_k$ and compute $\norm{r_k}_0$ using \eqref{e:residu0}. Note that one can reduce offline and online computational burden, as well as storage requirement, by noticing that $\Gamma_{j,j'}^c=\Gamma_{j',j}^c$ for all $j,j'$.

The cost of computation (and storage) of $\rho_k$ and $\norm{\rho_k}$ asymptotically dominates the cost of the online phase. This cost is in $O\left( \left( n_S(f)n_T(f) + N^2 \right)^2 \right)$. Again, we see that our weak Dirichlet treatment, as $n_T(f)$ remains independent of $n(b_0)$ and $n(b_1)$, allows for a better complexity of the online phase.

\subsection{Lower and upper bounds on stability constant}
\label{ss:scm}
To find lower and upper bounds on $C_k$ efficiently in order to use them in our computable error bound of Section \ref{sss:compt}, we turn to the successive constraints method (SCM) \cite{huynh2007successive,nguyen2009reduced}. Here we present the application to our case for the sake of self-containedness. Our difference is the use of the metric \eqref{e:metricSCM} during the constraint-selection phase.

\textbf{Notation. } As in Section \ref{ss:choicebasis}, we will need to handle reduced solutions of several values of the parameter tuple $\mu\in\mathcal P$ (see \eqref{e:definitionP}), for different timesteps $k=1,\ldots,\mathcal T$. We thus define an application that is the "reduced" counterpart of $u$ defined in Section \ref{ss:choicebasis}:
\[
\widetilde{u} : \left\{ \begin{aligned} \{1,\ldots,\mathcal T\} \times \mathcal P &\rightarrow \widetilde X_0 \\
\left( k, \mu \right) &\mapsto \widetilde{u}(k,\mu) = \utk \text{ satisfying \eqref{e:discrb2} for } \mu \text{ as parameter} \end{aligned} \right.
\]
We make $C_k$ depend explicitly on $\mu$ by defining:
\[ C_k(\mu) = \inf_{v\in X_0, \norm{v}=1} \left[   2 c (\widetilde u (k,\mu), v, v) + \nu(\mu) a (v,v)\right] \]
\textbf{SCM lower bound. } We now proceed to the derivation of the SCM lower bound of $C_k(\mu)$. We use the RB expansion:
$ \widetilde u(k,\mu) = \sum_{j=1}^N u_j^k(\mu) \zeta_j $ 
to rewrite $C_k(\mu)$ as:
\begin{align*}
C_k(\mu) = \inf_{v\in X_0, \norm v=1} \left[ \sum_{j=1}^N 2 u_j^k(\mu) c(\zeta_j,v,v) + \nu a(v,v) \right] 
&= \inf_{y=(y_1,\ldots,y_{N+1})\in\mathcal Y} \left[ \sum_{j=1}^N 2 u_j^k(\mu) y_j + \nu y_{N+1} \right] \\
&= \inf_{y \in \mathcal Y}  \mathcal J(\mu,k,y) 
\end{align*}
where:
\begin{multline*}
\mathcal Y = \big\{ y=(y_1,\ldots,y_{N+1}) \in \R^{N+1} \vert 
\exists v \in X_0, \norm v=1 \text{ s.t. } 
y_j=c(\zeta_j,v,v) \forall j=1,\ldots,N, y_{N+1}=a(v,v)
\big\} \end{multline*}
and:
\[ 
\mathcal J(\mu,k,y) = 2 \sum_{j=1}^N u_j^k(\mu) y_j + \nu y_{N+1} 
\]
We define, for a given "constraint subset" $\mathcal C \subset \{1,\ldots,\mathcal T\}\times P$ that will be chosen later: 
\begin{itemize}
\item 
$ \widetilde {\mathcal Y} = \left\{ (y_1,\ldots,y_{N+1}) \in \prod_{i=1}^{N+1} \left[ \sigma_i^{min}; \sigma_i^{max} \right] \vert
 \mathcal J(\mu',k',y) \geq C_{k'}(\mu'), \; \forall (\mu', k') \in \mathcal S(\mu,k) \right\} $
with $\mathcal S(\mu,k)$ standing for the set of the $M$ points in $\mathcal C$ that are closest to $(\mu,k)$ with respect to this metric:
\begin{equation}
\label{e:metricSCM}
d\big( (\mu,k), (\mu',k') \big) = \sum_{i=1}^{\dim\mathcal P} \left( \frac{ \mu^i - {\mu'}^i }{ \mu^i_{min} - \mu^i_{max} } \right)^2 + \left( \frac{ k - k'}{\mathcal T} \right)^2 
\end{equation}
where $\left( \mu^1, \ldots, \mu^{\dim\mathcal P} \right)$ are the coordinates of $\mu\in\mathcal P$, and:
\[ \mu^i_{min} = \min_{\mu\in\mathcal P} \mu^i, \;\;\;\;\; \mu^i_{max} = \max_{\mu\in\mathcal P} \mu^i \]
for $i=1,\ldots,\dim\mathcal P$ (here $\dim\mathcal P = 1+2+n(b_0)+n(b_1)+n_T(f) n_S(f)+n(u_0)$ is the number of parameters).

The metric defined in \eqref{e:metricSCM} quantifies proximity of two parameter-time tuples, with appropriate weighting so as to account for scaling differences between the parameters.

\item We further define:
\begin{align*}
\sigma_i^{min} &= \inf_{v\in X_0, \norm v=1} c(\zeta_i,v,v),\quad \forall i=1,\ldots,N &
\sigma_{N+1}^{min} &= \inf_{v\in X_0, \norm v=1} a(v,v) \\
\sigma_i^{max} &= \sup_{v\in X_0, \norm v=1} c(\zeta_i,v,v),\quad \forall i=1,\ldots,N & 
\sigma_{N+1}^{max} &= \sup_{v\in X_0, \norm v=1} a(v,v) 
\end{align*}
\end{itemize}
The SCM lower bound is then given by the following lemma:
\begin{lmm}[Proposition 1 in \cite{huynh2007successive}]
\label{l:1}
For every $\mathcal C \subset \{1,\ldots,\mathcal T\}\times P$ and $M\in\N$, and every $k=1,\ldots,\mathcal T$, a lower bound for $C_k(\mu)$ is given by:
\begin{equation}
\label{e:scmlow}
C_k^{inf}(\mu) = \inf_{y \in \widetilde{\mathcal Y}} \left[ \mathcal J(\mu,k,y) \right] 
\end{equation}
\end{lmm}

An algorithm for choosing a constraint subset $\mathcal C$ will be given after the description of the SCM upper bound and the SCM offline/online procedure.

\textbf{SCM upper bound. } We define:
\[ \widetilde {\mathcal Y}^{up} = \{ y^*(k_i, \mu_i) \; ; \; i=1,\ldots,I \; ; \; (k_i,\mu_i)\in\mathcal C \} \]
where:
\[ \mathcal C = \left\{ (k_1,\mu_1), (k_2,\mu_2), \ldots, (k_I,\mu_I) \right\} \]
and:
\[ y^*(k_i,\mu_i) = \arginf_{y\in\mathcal Y} \left[ \mathcal J(\mu_i, k_i, y) \right] \;\;\; (i=1,\ldots,I) \]
\begin{lmm}[Proposition 1 in \cite{huynh2007successive}]
\label{l:2}
For every $k=1,\ldots,\mathcal T$, an upper bound for $C_k$ is given by:
\begin{equation}
\label{e:scmup}
C_k^{sup}(\mu) = \inf_{y \in \widetilde{ \mathcal Y }^{up} }  \mathcal J(\mu,k,y)  
\end{equation}
\end{lmm}

\textbf{SCM offline/online procedure. } Relying on Lemmas \ref{l:1} and \ref{l:2}, our offline/online procedure for computing $C_k^{inf}$ and $C_k^{sup}$ reads:
\begin{algo}
\caption{SCM offline/online}
	\label{algo:3}
\begin{itemize}
\item offline:
 \begin{enumerate}
  \item choose $M$ and constraint set $\mathcal C$ (see next paragraph) ;
  \item compute and store $\sigma_i^{min}$ and $\sigma_i^{max}$ ($i=1,\ldots,N+1$) by solving a generalized eigenproblem on $X_0$ ;
  \item for each $(k',\mu')\in\mathcal C$:
  	\begin{enumerate}
	\item solve a generalized eigenproblem on $X_0$ to find $C_{k'}(\mu')$ (and store it);
	\item let $w\in X_0$ be a unit eigenvector of the above eigenproblem;~compute and store (in $\widetilde{\mathcal Y}^{up}$) $y^*(k',\mu')$ using:
	 \begin{align*}
			y^*(k',\mu')_j &= c(\zeta_j,w,w) \;\; (j=1,\ldots,N) \\
			y^*(k',\mu')_{N+1} &= a(w,w)
	 \end{align*}
	 \end{enumerate}
  \end{enumerate}
\item online: 
  \begin{itemize}
  \item for the lower bound:
  \begin{enumerate}
   \item assemble and solve optimization problem \eqref{e:scmlow} ;
	\end{enumerate}
	\item for the upper bound: test one-by-one each element of $\widetilde{ \mathcal Y }^{up} $ to solve \eqref{e:scmup}.
   \end{itemize}
\end{itemize}
\end{algo}
In the lower bound online phase, the optimization problem required to solve is a \emph{linear programming} problem (LP) with $N+1$ variables and $N+1+M$ constraints ($M$ one-sided inequalities and $N+1$ two-sided). There are algorithms, such as the simplex algorithm (see \cite{nocedal1999numerical} for instance), which solve such optimization problems under (on average) polynomial complexity with respect to the number of variables and number of constraints, even if they can be exponential in the worst cases. What matters here is this complexity is independent of $\mathcal N$. The upper bound online phase has a complexity depending linearily on the cardinality of the reasonably-sized $\mathcal C$ and on $N$.

\textbf{"Greedy" constraint set selection. }
To choose $\mathcal C$ in Algorithm \ref{algo:3}, step 1, we can use the greedy constraint set selection Algorithm \ref{algo:4}.
\begin{algo}
\caption{Greedy constraint set selection}
	\label{algo:4}
\begin{enumerate}
\item choose $M\in\N$ ;
\item initialize $\mathcal C = \{ (k_1, \mu_1) \}$ with arbitrary $k_1 \in \{1,\ldots,\mathcal T\}$ and $\mu_1\in\mathcal P$;
\item choose a rather large, finite-sized sample $\Xi \subset \{1,\ldots,\mathcal T\} \times \mathcal P$ ;
\item repeat:
	\begin{itemize}
	\item using the "current" $\mathcal C$ to compute $C^{sup}$ and $C^{inf}$, append:
	\[ (k^*, \mu^*) = \argmax_{(k,\mu)\in\Xi} \frac{ \exp\left(C_k^{sup}(\mu)\right) - \exp\left(C_k^{inf}(\mu)\right) }{\exp\left(C_k^{sup}(\mu)\right)} \quad\quad\text{ to } \mathcal C. \]
	\end{itemize}
\end{enumerate}
\end{algo}
The repeat loop in this algorithm can be stopped either when $\#\mathcal C$ has reached a maximal value, or when the "relative exponential sharpness" indicator:
\[ \max_{(k,\mu)\in\Xi} \frac{ \exp\left(C_k^{sup}(\mu)\right) - \exp\left(C_k^{inf}(\mu)\right) }{\exp\left(C_k^{sup}(\mu)\right)} \]
gets less than a desired precision. We use a measure of the difference between the exponentials so as to account for the "exponential" effect of the stability constants on the error bounds \cite{nguyen2009reduced}.

As in the greedy algorithm for basis selection described at Section \ref{sss:greedy}, this algorithm makes, at each step, the "best possible" choice, that is the value of the parameter and time for which the bounds computed using the current constraint set are the less sharp.

A last remark we can give on the algorithm is about the trade-off in the choice of $M$: whatever $M$ is, we always get a certified bound on $C_k(\mu)$, but increasing $M$ will improve sharpness of this bound, at the expense of an increase in online computation time.

\section{Numerical results}
\label{s:numres}
We now present some numerical results obtained with the methodology described above. We implemented it in a software package written in C++, using GNU OpenMP \cite{gomp} as threading library, ARPACK \cite{arpack96} for eigenvalues computation and GLPK \cite{glpk10} as linear programming problems solver.

For all the experiments above, the convergence test for Newton iterations when solving \eqref{e:nleq} and \eqref{e:discrb2} was the following: $\norm{\delta}^2 \leq 3 \times 10^{-16}$. The penalization constant used was $P=10^7$.

We also took $M = \#\mathcal C = 10$ as parameters for the SCM procedure.

The $\Phi$ functions appearing in the parametrizations of $u_0$, $b_0$, $b_1$ and $f$ are chosen to be sine functions with fixed known angular velocity. More specifically, we suppose that:
\begin{align*}
b_0(t) &= b_{0m} + \sum_{l=1}^{n(b_0)} A^{b_0}_l \sin( \omega^{b_0}_l t ) & 
b_1(t) &= b_{1m} + \sum_{l=1}^{n(b_1)} A^{b_1}_l \sin( \omega^{b_1}_l t ) \\
f(t,x) &= f_m + \sum_{l=1}^{n_T(f)} \sum_{p=1}^{n_S(f)} A^f_{lp} \sin( \omega^{fT}_l t ) \sin( \omega^{fS}_p x) &
u_0(x) &= u_{0m} + \sum_{l=1}^{n(u_0)} A^{u_0}_l \sin( \omega^{u_0}_l x )
\end{align*}

\subsection{Reference solutions}
\begin{figure} 
\begin{center}
\includegraphics[scale=.7]{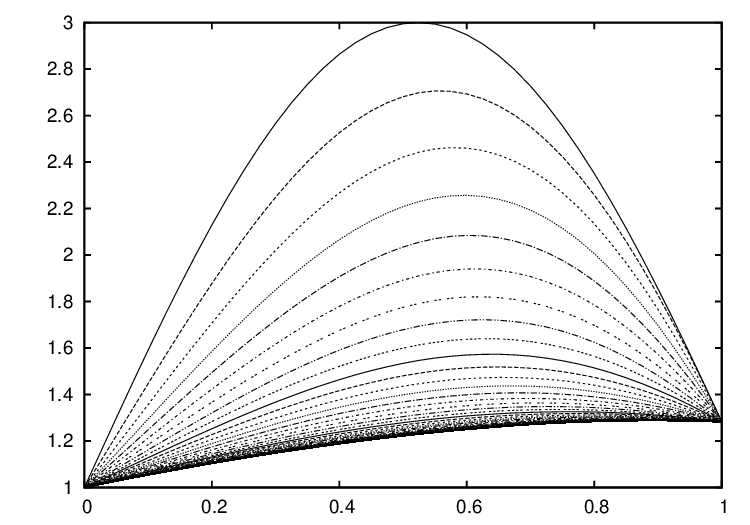}
\caption{Full solution with high viscosity $\nu=1$. Plots of the solution $u$ of equation \eqref{e:burg}, \eqref{e:init}, \eqref{e:bound2} as a function of space $x$, for various times $t$ ranging from $t=0$ to $t=2$ (the bottom lines correspond to high times). We use the parameters defined in \eqref{e:params1}. }
\label{f:1}
\end{center}
\end{figure}

Figure \ref{f:1} shows an example of the reference solution of \eqref{e:burg}, \eqref{e:init} and \eqref{e:bound2} every 10 timesteps. The parameters were:
\begin{equation}
\label{e:params1}
\begin{aligned}
\mathcal N &= 40 & 
\delt &= .02 \\
T &= 2 &
\nu &= 1 \\
b_0(t) &= 1 &
b_1(t) &= 1.28224 \\
f(t,x) &= 1 &
u_0(x) &= 1 + 2 \sin(3x)
\end{aligned}
\end{equation}

\begin{figure} 
\begin{center}
\includegraphics[scale=.7]{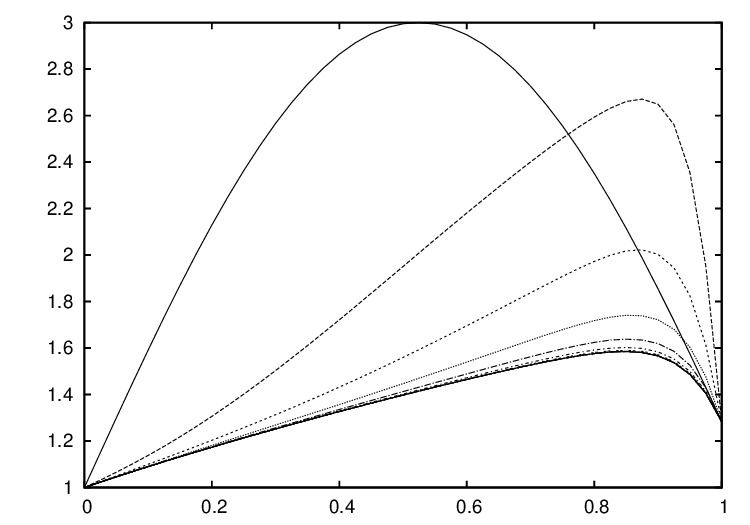}
\caption{Full solution with low viscosity $\nu=.1$. Plots of the solution $u$ of \eqref{e:burg}, \eqref{e:init}, \eqref{e:bound2} as a function of space $x$, for various times $t$ ranging from $t=0$ to $t=2$ (the bottom lines correspond to high times). We use the parameters defined in \eqref{e:params2}. }
\label{f:2}
\end{center}
\end{figure}

Figure \ref{f:2} does the same with a lower viscosity. The parameters were:
\begin{equation}
\label{e:params2}
\begin{aligned}
\mathcal N &= 40 &  
\delt &= .002 \\
T &= 2 &
\nu &= .1 \\
b_0(t) &= 1 &
b_1(t) &= 1.28224 \\
f(t,x) &= 1 &
u_0(x) &= 1 + 2 \sin(3x) 
\end{aligned}
\end{equation}

Figures \ref{f:1} and \ref{f:2} show the solution $u$ of the viscous Burgers' equation plotted as functions of space $x$, for various times $t$, respectively for the parameter set \eqref{e:params1} and \eqref{e:params2}. 

We have checked that the \emph{a posteriori} indicator of boundary error \eqref{e:posterioridiagnosis} gets no higher than $6 \times 10^{-7}$ in both cases. 

\subsection{Reduced solutions}
\subsubsection*{Computational economy}
\label{ss:economy}
To show the substantial time savings provided by the reduced basis approximation, we compute the reduced solution for the parameters  set given by \eqref{e:params1}, with $\mathcal N=60$. The full solution (with $\mathcal N=60$) takes 0.26s CPU time to be produced (when using Thomas' algorithm for tridiagonal matrices inversion).

\begin{figure} 
\begin{center}
\includegraphics[width=\textwidth]{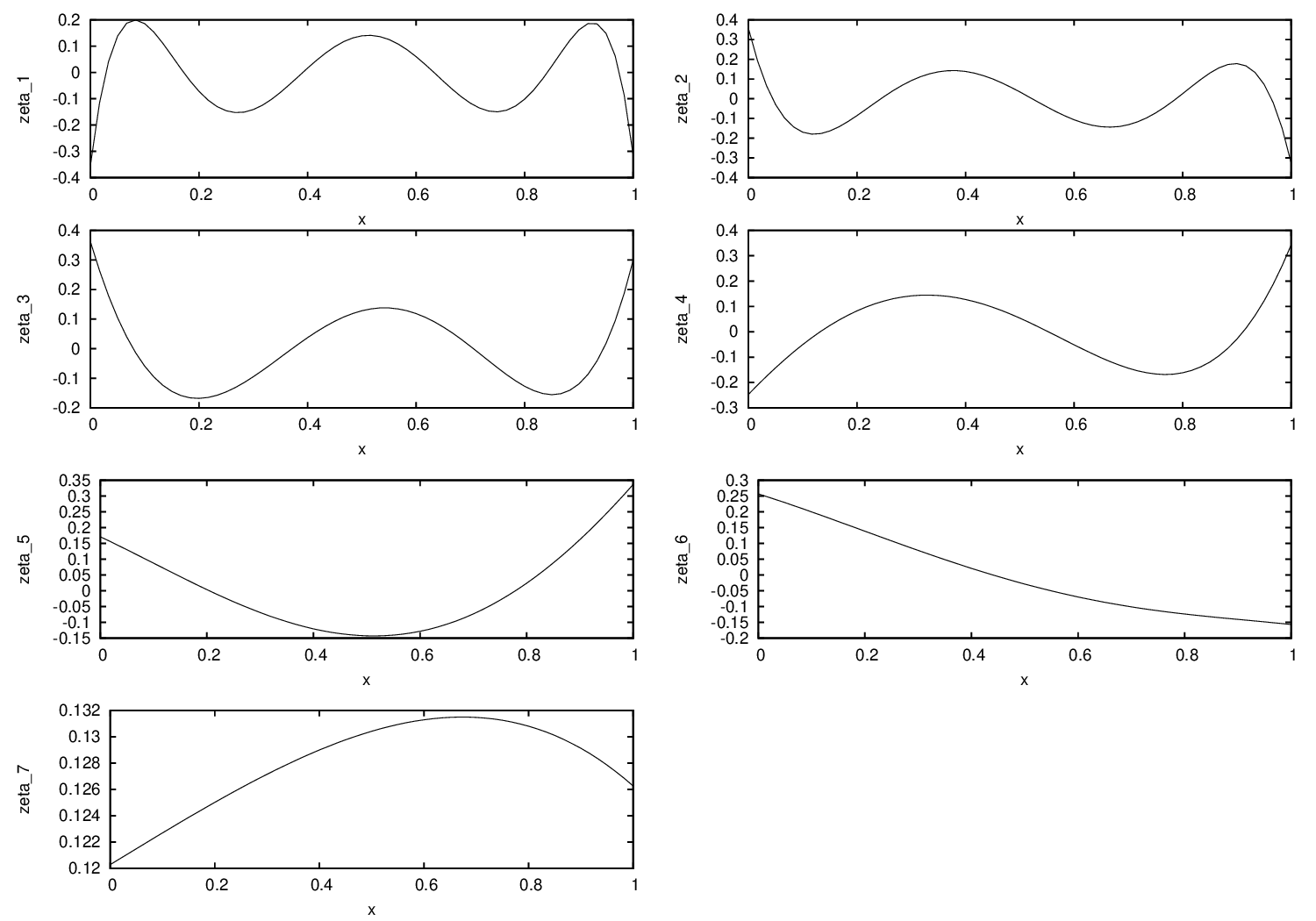}
\caption{POD reduced basis: plots, as functions of space, of the 7 leading POD modes, i.e. the $\zeta_i$ ($i=1,\ldots,7$) defined by \eqref{e:podeqn}, with $z_i$ ($i=1,\ldots,7$) the leading eigenvectors of $M^T\Omega M$. The modes are sorted (from top to bottom, and from left to right) by increasing magnitude of eigenvalues. Parameter ranges for snapshot sampling are those in Table \ref{tbl:1}. }
\label{f:3}
\end{center}
\end{figure}

We use the POD-driven basis selection procedure to select the $N=7$ leading POD modes (using $S=30$ snapshots). The resulting basis  (with the functions sorted by increasing magnitude of eigenvalues) is shown in Figure \ref{f:3}. We did not make use of the "expansion" procedure described in Section \ref{sss:enrich}. The overall CPU time for the offline phase was 6.36s.

We used fixed parameters $n(b_0)=n(b_1)=n_S(f)=n_T(f)=n(u_0)=1$ and parameter ranges as shown in Table \ref{tbl:1}.
\begin{table}[h]
\begin{center}
\begin{tabular}{|r|r|r||r|r|r|}
\hline
Parameter & Min. & Max. & Parameter & Min. & Max. \\
\hline
$\nu$& .8 & 1.2 & $A_1^{u_0}$ & 1.1 & 3 \\
$A_1^{b_0}$ & .9 & 1.2 &$\omega_1^{b_0}$ & 1 & 1 \\
$A_1^{b_1}$ & .9 & 1.2 &$\omega_1^{b_1}$ & 1 & 1 \\
$f_m$ & 0 & 2 & $\omega_1^{fT}$ & 2 & 2 \\
$A_{1,1}^f$ & 0.7 & 1.3 & $\omega_1^{fS}$ & 2 & 2 \\
$u_{0m}$ & 0 & 1 & $\omega_1^{u_0}$ & 3 & 3 \\
\hline
\end{tabular}
\end{center}
\caption{\label{tbl:1} Ranges of the different parameters in the benchmark of Section \ref{ss:economy}.}
\end{table}

\begin{figure} 
\begin{center}
\includegraphics[scale=.7]{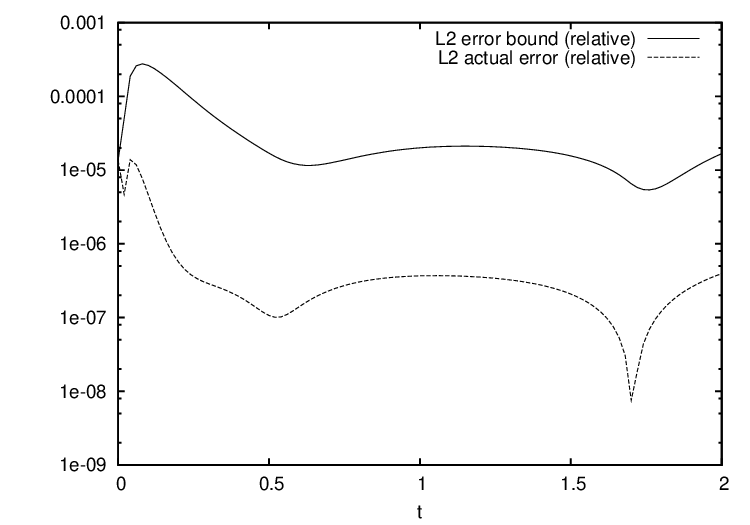}
\caption{Relative $L^2$ online error bound and actual error. We plot in solid line $\frac{\epsilon_k}{\norm{\utk}}$, and in dashed line $\frac{\norm{u^k-\utk}}{\norm{\utk}}$ as functions of $t=k\delt$ for $k=1,\ldots,\mathcal T$. }
\label{f:a1}
\end{center}
\end{figure}

We then used this basis to compute the reduced solution for a particular (randomly chosen) instance of the parameters. The reduced solution was computed in 0.04s, \emph{including} the time necessary for the online error bound computation, shown in Figure \ref{f:a1} (solid line). Our procedure reduces the marginal cost to 15\% of the original cost, yet providing a certified $L^2$ relative error of less than 1\textperthousand. 

\subsubsection*{Error bound estimation}
\begin{figure} 
\begin{center}
\includegraphics[scale=.7]{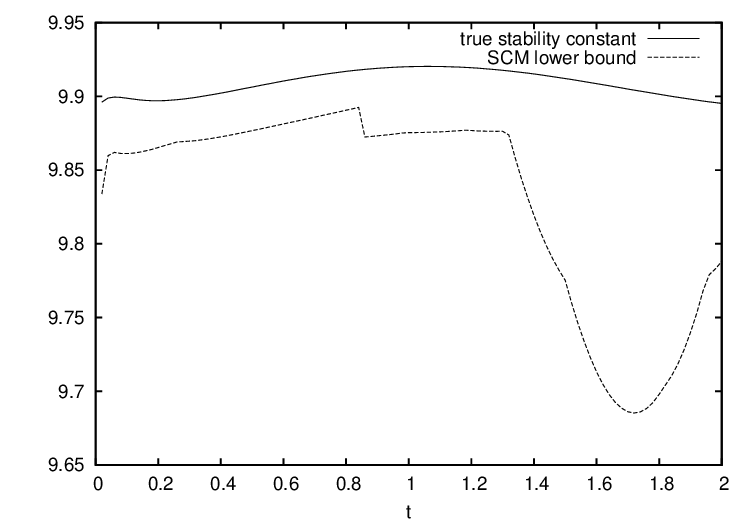}
\caption{True stability constant $C_k$ defined by \eqref{e:defCk} and SCM lower bound $C_k^{inf}$ defined by \eqref{e:scmlow} as functions of time. We use $M=\#\mathcal C=10$ as SCM parameters. }
\label{f:a2}
\end{center}
\end{figure}

Still using the preceding POD basis and instance of the parameters, we compared the online error bound with the actual error, for the same parameter set as above. The result is shown in Figure \ref{f:a1}. We see that our error bound is quite sharp, especially when it follows the decrease in the actual error near $t=1.3$. We also checked for the quality of the SCM procedure, by comparing the actual stability constant with the lower bound provided by SCM (Figure \ref{f:a2}).

\subsection{Convergence benchmarks}
In order to compare our three basis selection procedures (POD, greedy and POD-Greedy), we have made "convergence benchmarks", i.e. representations of the maximal and mean (estimated) error over all timesteps, and over a sample of 100 parameters as functions of the size of the reduced basis. The same sample of benchmark parameters is used throughout all the procedure.

\begin{figure} 
\begin{center}
\includegraphics[scale=.7]{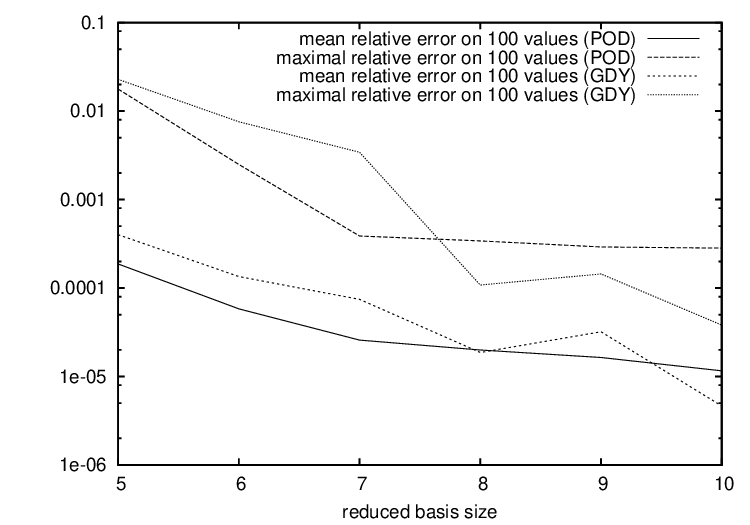}
\caption{Convergence benchmark 1. We plot (on a logarithmic scale) maximal and mean relative online error bounds over a (uniform) random sample of 100 initial values $u_{0m}\in[0,1]$, as functions of the reduced basis size $N$, when the reduced basis is chosen using POD-based procedure (POD) or greedy procedure (GDY). Fixed parameters are $\nu=1$ and $f_m=1$. }
\label{f:5}
\end{center}
\end{figure}

Comparison of greedy (with $\#\Xi=100$) and POD (with $S=60$) procedures for $\mathcal N=40$, $T=2$, $\delt=.02$, $n(b_0)=n(b_1)=n_S(f)=n_T(f)=n(u_0)=0$, with parameters $f_m$ and $\nu$ fixed to unity, and varying $u_{0m}\in [0,1]$ (and thus initial boundary values $b_{0m}$ and $b_{1m}$, moving accordingly to compatibility conditions \eqref{e:compatcond}) is shown in Figure \ref{f:5}. 

The benchmarking process for greedy took 19.5s of CPU time, the one for POD took 14.98s. The online cost, depending only on the size of the reduced basis, is the same regardless of how the reduced basis has been chosen. We also see the fast (exponential) convergence of error bound towards zero as $N$ increases.

\begin{figure} 
\begin{center}
\includegraphics[scale=.7]{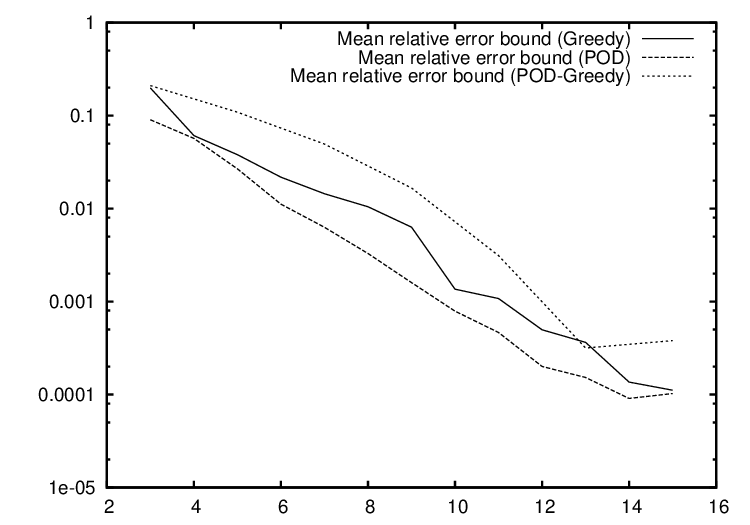}

\includegraphics[scale=.7]{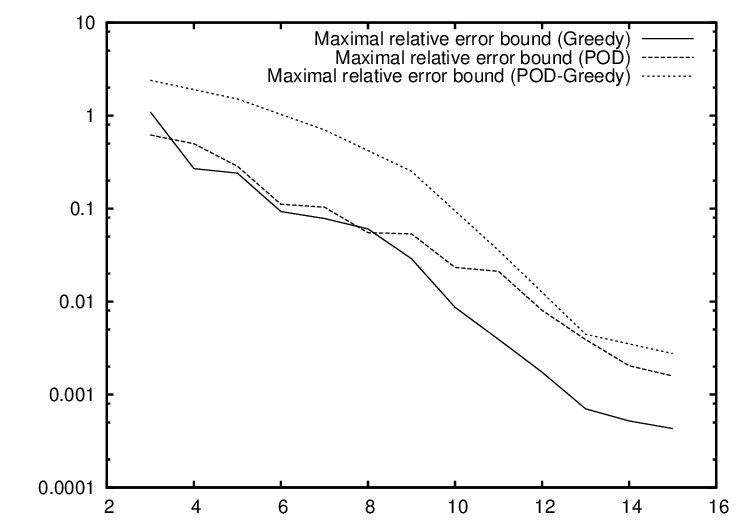}
\caption{Convergence benchmark 2. We plot (on a logarithmic scale) maximal (bottom) and mean (top) online error bounds over a (uniform) random sample of 100 initial values $u_{0m}\in[0,1]$, as functions of the reduced basis size $N$, when reduced basis is chosen using POD-based procedure with $S=90$, Greedy, or POD-Greedy procedure with $P_1=2$. Fixed parameters are $\nu=0.1$ and $f_m=1$. }
\label{f:7}
\end{center}
\end{figure}

\begin{figure} 
\begin{center}
\includegraphics[width=\textwidth,height=.7\textheight]{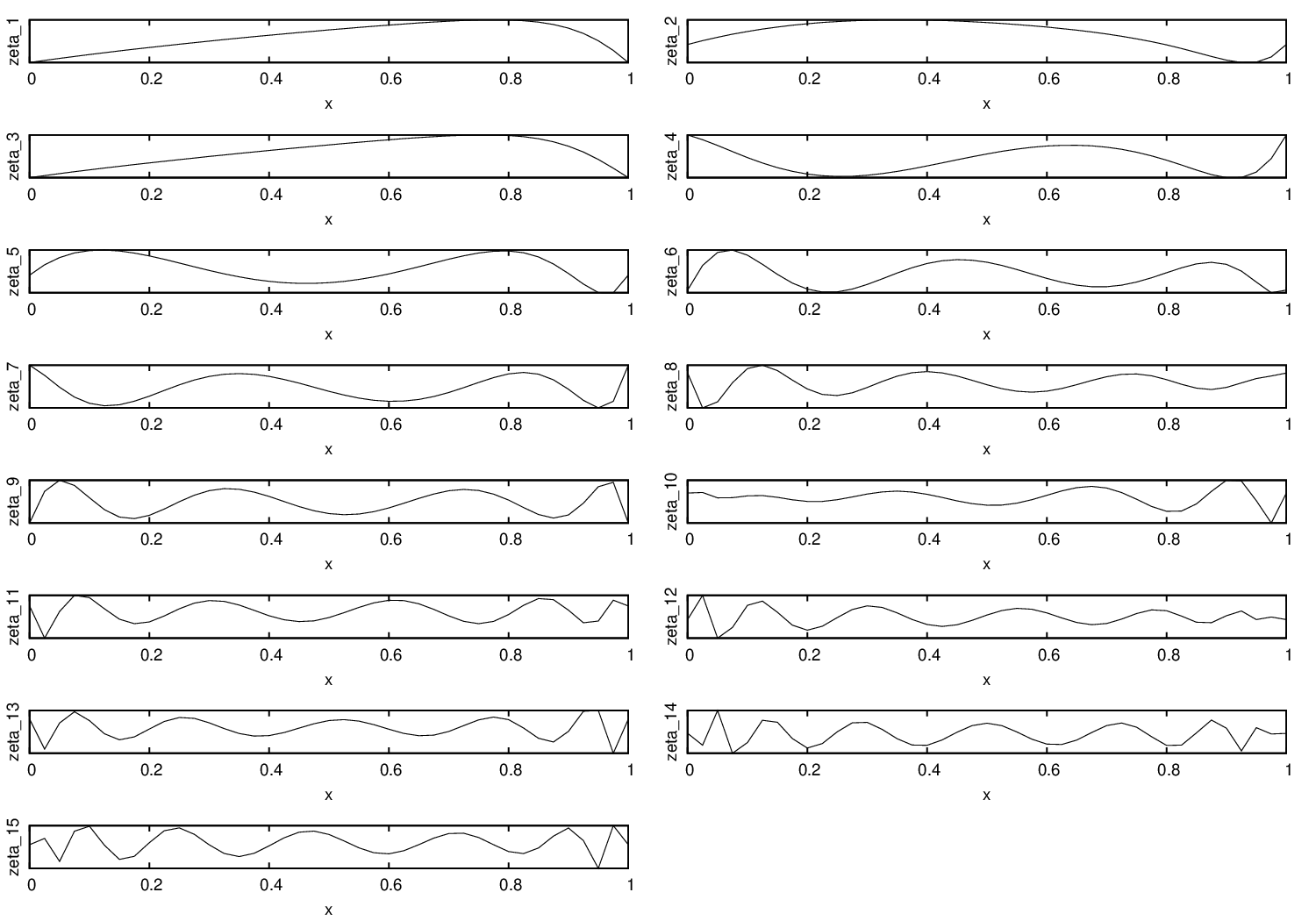}
\caption{Reduced basis selected at the final step of the POD-Greedy procedure carried out for the previous benchmark (Figure \ref{f:7}). Basis elements are plotted as functions of space, and are sorted (from top to bottom, left to right) in order of selection in the greedy procedure. }
\label{f:6}
\end{center}
\end{figure}

Another benchmark was then made, with the same data, except that $\nu=.1$ and $\delt=.002$.  
The result is visualized in Figure \ref{f:7}. The POD benchmarking process took 349 s of CPU time, the POD-Greedy, 282 s, and the Greedy, 470 s. We notice that a smaller viscosity leads to degraded precision of our RB approximation. 
The resulting final basis selected by the POD-Greedy is displayed in Figure \ref{f:6}. 

The POD-Greedy algorithm is run using $P_1=2$ and initialized using the full time-discrete trajectory (hence, 200 vectors) for one parameter value. Using the initial data as initialization for the POD-Greedy algorithm and using $P_1=1$ may give better overall performance of POD-Greedy (at the expense of an increased offline computation time), as the sample of error bounds is updated at every step of the algorithm (instead of every other step for $P_1=2$).

\subsection{Effect of mesh refinement and penalization constant}
To check for the robustness of our bound, we studied the influence of $\mathcal N$ (the number of ``full'' spatial discretization points) and $P$ on the magnitude and the sharpness of our error bound. We ran the same benchmark as in Section \ref{ss:economy}, but with different $\mathcal N$ or $P$.

We did the test with refined meshes ($\mathcal N=200$, $\mathcal N=800$) and we obtained a similar error bound profile; we conclude that the sharpness of our error bound is quite insensitive to mesh refinement.

\begin{figure} 
\begin{center}
\includegraphics[scale=.6]{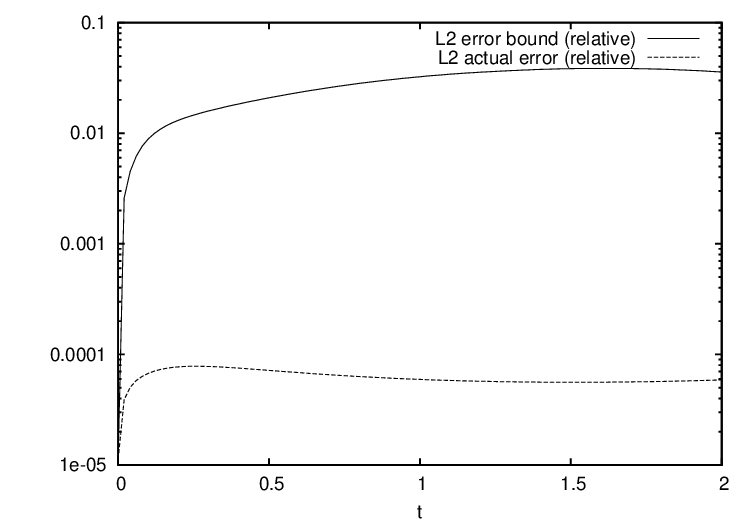}
\includegraphics[scale=.6]{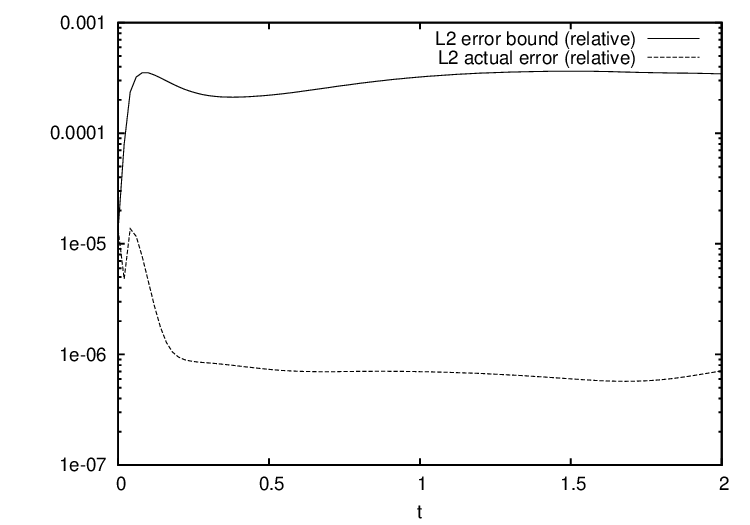}

\includegraphics[scale=.6]{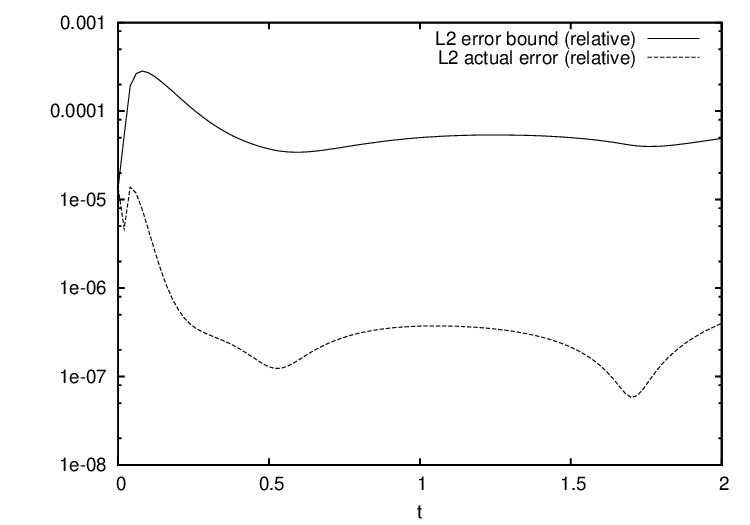}
\includegraphics[scale=.6]{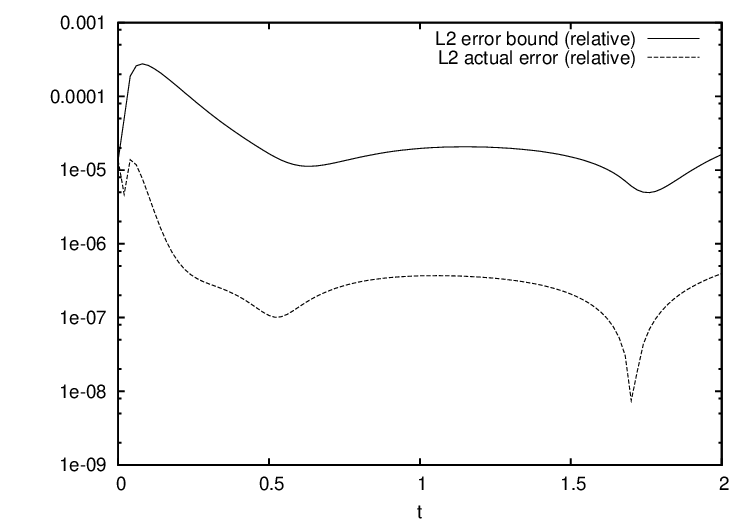}
\caption{Relative $L^2$ online error bound and actual error. We plot in solid line $\frac{\epsilon_k}{\norm{\utk}}$, and in dashed line $\frac{\norm{u^k-\utk}}{\norm{\utk}}$ as functions of $t=k\delt$ for $k=1,\ldots,\mathcal T$, and $\mathcal N=60$ and, from left to right and top to bottom, $P=10^2$, $P=10^4$, $P=10^5$ and $P=10^{12}$.}
\label{f:10}
\end{center}
\end{figure}

In Figure \ref{f:10}, we visualize the actual error and the error bound for various values of $P$ (with $\mathcal N=60$). We see that the error bound is tighter and sharper for high values of the penalization constant $P$. This can easily be explained by the fact that, as $P \rightarrow +\infty$, the errors at the boundary $e_k(0)$ and $e_k(1)$ vanish, hence modifying all the terms in the error bound (i.e., decreasing $\eta_k$, $\sigma_k$, $f_k$, $\xi_k^{\mathcal A}$, $\xi_k^{\mathcal B}$, $\xi_k^{\mathcal \gamma}$, $\mathcal B_k$, $\gamma_k$  and increasing $\mathcal A_k$) where these errors appear.


\subsection{Comparison with existing bound}
\begin{figure} 
\begin{center}
\includegraphics[scale=.7]{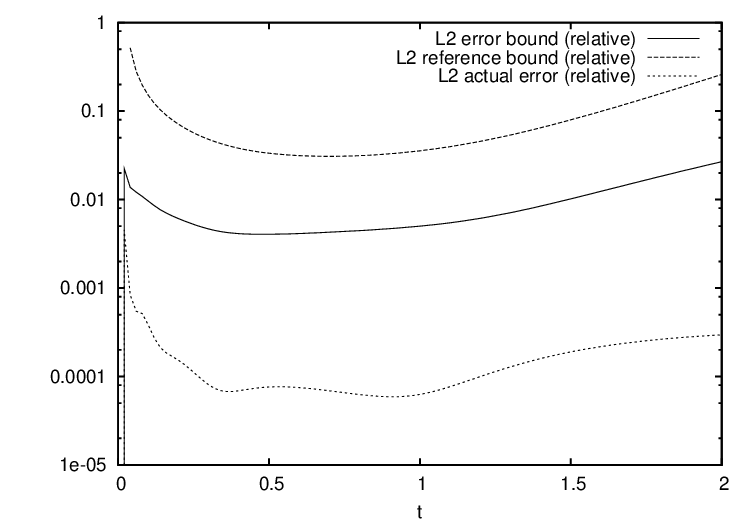}
\caption{Comparison with the existing error bound for reduced basis Burgers equation (reference bound). We plot the actual error, the existing (reference) error bound and our error bound as functions of time. We took $\nu=0.1$. }
\label{f:compa}
\end{center}
\end{figure}

We compared our error bound with the bound described in \cite{nguyen2009reduced}, for $\mathcal N=60$, $\Delta t=.02$, $ T=2$, fixed $u_0=b_0=b_1=0$, $f=1$ and variable $\nu\in[0.1;1]$. The reduced basis is found by POD with $S=90$. Figure \ref{f:compa} shows that our bound is clearly better than the existing reference bound.
 
Besides, we have performed a benchmark over a fixed sample of 100 random values of $\nu$  uniformly chosen in $[0.1;1]$. For our bound, the mean error bound is $0.00076$, the maximum error bound is $0.02$, while for the bound of \cite{nguyen2009reduced}, we obtain $0.0041$ for the mean bound and $0.25$ for the maximum.

It is easy to see that, when the boundary conditions are fixed to zero, the recurrence formula for the error bound reduces to:
	\[ \norm{e_k} \leq \frac{ \norm{e_{k-1}}+\Delta t \norm{r_k}_0 }{1+C_k \Delta t} \]
	while the error bound described in \cite{nguyen2009reduced} has the following expression:
	\[ \norm{e_k} \leq \sqrt{ \frac{ \norm{e_{k-1}}^2+\frac{\Delta t}{\nu} \norm{r_k}_0^2 }{1+ \widetilde{C_k} \Delta t } } \]
	where the modified stability constant $\widetilde{C_k}$ reads:
	\[ \widetilde{C_k} = \inf_{v\in X_0} \frac{4 c(\widetilde u^k,v,v)+\nu a(v,v)}{\norm{v}^2} \]
	We recall that our stability constant $C_k$ is given by:
	\[ C_k= \inf_{v\in X_0, \norm{v}=1} \left[ 2c(\utk,v,v)+\nu a(v,v) \right] \]
	The $\nu$-dependence of the bound can explain why our bound is better, especially for small values of $\nu$. Besides, the derivation of our error bound makes lesser use of inequalities (for instance, we do not make use of Young's inequality at the beginning of the proof, each inequality used is a potential source of optimality loss) and keeps treating more terms.

\section*{Conclusion}
We have presented a certified procedure for low marginal cost approximate resolution of the viscous Burgers equation with parametrized viscosity, as well as initial and boundary value data. This procedure makes use of a reduced basis offline/online procedure for a penalized weak formulation, an efficiently computed error bound in natural $L^2$ norm (made possible by the successive constraints method (SCM)), and three procedures at hand for choosing a basis to expand reduced solutions in.

Our procedure becomes less useful when the ratio time/viscosity increases, as this degradates the stability constant $C_k$. Another limitation of our method, for one willing to use it for large times, is that the online procedure complexity still depends on the temporal discretization step. However, our numerical experiments show a substantial decrease in marginal cost when using reduced basis approximation, as well as efficiency (both in terms of sharpness and computation time) of the provided error bound for moderate viscosities. This decrease in the cost is made possible by the fact that online procedure has a complexity that is independent from the number of spatial discretization points. 

\appendix
\section{Proof of Theorem \ref{thm:1}}
\begin{proof}
Subtracting left-hand side of \eqref{e:discrb2} from both sides of relation \eqref{e:uktruth} yields that for $k=1,\ldots,\mathcal T$, the error at time $t_k$: $e_k=u^k_e-\utk$ satisfies, for every $v\in X_0$:
\begin{equation}
\label{e:proof1}
\frac{1}{\Delta t} ( \langle e_k, v \rangle - \langle e_{k-1},v \rangle ) + c(u^k_e,u^k_e,v)-c(\utk,\utk,v) +
 \nu a(e_k,v) = r_k(v) 
\end{equation}
We write:
\begin{equation}
\label{e:ekdecomp}
e_k = e_k(0) \phi_0 + e_k(1) \phi_{\mathcal N} + e_k^z 
\end{equation}
with $e_k^z \in X_0 = \{ v \in X \text{ st. } v(0)=v(1)=0 \}$.
And, by applying \eqref{e:proof1} with $v=e_k^z$:
\begin{equation}
\label{e:etoile}
\frac{1}{\Delta t} ( \langle e_k,e_k^z \rangle - \langle e_{k-1},e_k^z \rangle ) + c(u^k_e,u^k_e,e_k^z)-c(\utk,\utk,e_k^z) +
 \nu a(e_k,e_k^z) = r_k(e_k^z) 
\end{equation}
Since:
\begin{align*}
	\langle e_k,e_k^z \rangle &=e_k(0) \langle \phi_0, e_k^z \rangle + e_k(1) \langle \phi_{\mathcal N},e_k^z \rangle +  \norm{ e_k^z }^2 \\
	&= \norm{e_k^z}^2 + e_k(0) e_k^z\left(\frac{1}{\mathcal N}\right) \langle \phi_0,\phi_1 \rangle + e_k(1) e_k^z\left(1-\frac{1}{\mathcal N}\right) \langle \phi_{\mathcal N-1},\phi_{\mathcal N} \rangle
\end{align*}
and that, for every $v\in X$:
\begin{align*}
c(u^k_e,u^k_e,v)-c(\utk,\utk,v) &= -\frac{1}{2}\int_0^1 \left( (u^k_e)^2 - (\utk)^2 \right) \frac{\partial v}{\partial x} \\
&= -\frac{1}{2}\int_0^1  \left( u^k_e - \utk \right) \left( u^k_e + \utk \right) \frac{\partial v}{\partial x} \\
&= -\frac{1}{2}\int_0^1 e_k \left( \utk + \utk + e_k \right) \frac{\partial v}{\partial x} \\
&= 2c(\utk,e_k,v)+c(e_k,e_k,v)
\end{align*}
We have that \eqref{e:etoile} implies:
\begin{multline}
\label{e:etoile2}
\frac{1}{\Delta t} ( \norm{e_k^z}^2 - \langle e_{k-1},e_k^z \rangle ) + \psi_k(e_k,e_k^z) + \frac{1}{\Delta t} \left( e_k(0) e_k^z\left(\frac{1}{\mathcal N}\right) \langle \phi_0,\phi_1 \rangle + e_k(1) e_k^z \left(1-\frac{1}{\mathcal N}\right) \langle \phi_{\mathcal N},\phi_{\mathcal N-1} \rangle \right)
\\ = r_k(e_k^z) - c(e_k,e_k,e_k^z) 
\end{multline}

We are now willing to find a lower bound for the left-hand side of \eqref{e:etoile2} and an upper bound for its right-hand side.

\emph{Lower bound for LHS. }
From Cauchy-Schwarz inequality:
\[ -\langle e_{k-1},e_k^z \rangle \;\;\; \geq \;\;\; - \norm{e_{k-1}} \norm{e_k^z} \]
and, by the triangle inequality:
\begin{equation}
\label{e:triangle}
 \norm{e_k} - \eta_k \leq \norm{e_k^z} \leq \norm{e_k} + \eta_k 
\end{equation}
because $\eta_k = \abs{e_k(0)} \norm{\phi_0} + \abs{e_k(1)} \norm{\phi_{\mathcal N}}$.

So:
\begin{equation} 
\label{e:xx1}
-\langle e_{k-1},e_k^z \rangle \;\;\; \geq \;\;\; - \norm{e_{k-1}} ( \norm{e_k} + \eta_k )
\end{equation}
We also have:
\begin{equation}\label{e:ineg1} \norm{e_k^z}^2 \geq \norm{e_k}^2+\norm{e_k-e_k^z}^2-2\norm{e_k}\norm{e_k-e_k^z} \geq \norm{e_k}^2-2 \norm{e_k} \eta_k \end{equation}
because of Cauchy-Schwarz and triangle inequalities. 

Besides,
\[ 0 \leq \norm{e_k^z} \leq \norm{e_k^z-e_k} + \norm{e_k} \leq \eta_k + \norm{e_k} \]
so that:
\begin{equation}\label{e:ajout} \norm{e_k^z}^2 \leq \left(\norm{e_k}+\eta_k\right)^2 \end{equation}
Using the bilinearity of $\psi_k$, we have:
\begin{equation}
\label{e:xx2}
\psi_k(e_k,e_k^z) = \psi_k(e_k^z,e_k^z) + e_k(0) \psi_k(\phi_0,\phi_1)e_k^z\left(\frac{1}{\mathcal N}\right)+e_k(1) \psi_k(\phi_{\mathcal N},\phi_{\mathcal N-1}) e_k^z\left(1-\frac{1}{\mathcal N}\right) 
\end{equation}
because $\psi_k(\phi_0,\phi_{\mathcal N})=\psi_k(\phi_{\mathcal N},\phi_0)=0$, since $\phi_0$ and $\phi_{\mathcal N}$ have no common support, $\psi_k(\phi_0,\phi_j)=0$ for $j>1$, and $\psi_k(\phi_{\mathcal N},\phi_j)=0$ for $j<\mathcal N-1$.

From the definition of the stability constant $C_k$:
\[ \psi_k(e_k^z,e_k^z) \geq C_k \norm{e_k^z}^2 \]
So that, thanks to \eqref{e:ineg1} and \eqref{e:ajout},
\[ \psi_k(e_k^z,e_k^z) \geq \begin{cases} C_k \norm{e_k}^2-2\eta_k C_k \norm{e_k} & \text{ if } C_k\geq0 \\
		C_k \norm{e_k}^2+2\eta_k C_k \norm{e_k} + C_k \eta_k^2 & \text{ if } C_k\leq0 \end{cases} \]
That is
\begin{equation}
\label{e:xx3}
\psi_k(e_k^z,e_k^z) \geq C_k \norm{e_k}^2 - \sigma_k \norm{e_k} - \partieneg{C_k} \eta_k^2 
\end{equation}
from the definition of $\sigma_k$.

We have:
\begin{equation}
\label{e:e1}
\left| e_k^z\left( \frac{1}{\mathcal N} \right) \right| \leq \mathcal E \norm{e_k^z} \leq \mathcal E \norm{e_k} + \mathcal E \eta_k 
\end{equation}
and by symmetry:
\begin{equation}
\label{e:e2}
\left| e_k^z\left( 1 - \frac{1}{\mathcal N} \right) \right| \leq \mathcal E \norm{e_k^z} \leq \mathcal E \norm{e_k} + \mathcal E \eta_k 
\end{equation}
so that, combining \eqref{e:e1}, \eqref{e:e2} and introducing $f_k$:
\[
\abs{ e_k(0) \psi_k(\phi_0,\phi_1)e_k^z\left(\frac{1}{\mathcal N}\right)+e_k(1) \psi_k(\phi_{\mathcal N},\phi_{\mathcal N-1}) e_k^z\left(1-\frac{1}{\mathcal N}\right)} \leq
 \norm{e_k} f_k + \eta_k f_k 
\]
Now we can say that:
\begin{equation}
\label{e:xx4}
e_k(0) \psi_k(\phi_0,\phi_1)e_k^z\left(\frac{1}{\mathcal N}\right)+e_k(1) \psi_k(\phi_{\mathcal N},\phi_{\mathcal N-1}) e_k^z\left(1-\frac{1}{\mathcal N}\right) \geq - \norm{e_k} f_k - \eta_k f_k 
\end{equation}
And we also have:
\begin{align}
\label{e:xx5}
\abs{ e_k(0) e_k^z\left(\frac{1}{\mathcal N}\right)\langle \phi_0,\phi_1 \rangle+e_k(1) e_k^z\left(1-\frac{1}{\mathcal N}\right)  \langle \phi_{\mathcal N-1},\phi_{\mathcal N} \rangle}& \leq \abs{e_k(0)} \left(\mathcal E \norm{e_k}+\mathcal E \eta_k \right) \langle \phi_0,\phi_1 \rangle \notag\\
&\;\;\;\; + \abs{e_k(1)}\left(\mathcal E \norm{e_k}+\mathcal E \eta_k\right) \langle \phi_{\mathcal N},\phi_{\mathcal N-1} \rangle\notag\\
&=\mathcal E(\abs{e_k(0)}\langle \phi_0,\phi_1 \rangle+\abs{e_k(1)}\langle \phi_{\mathcal N},\phi_{\mathcal N-1} \rangle)\norm{e_k} \notag\\
&\;\;\;\; +\mathcal E \eta_k (\abs{e_k(0)}\langle \phi_0,\phi_1 \rangle+\abs{e_k(1)}\langle \phi_{\mathcal N},\phi_{\mathcal N-1} \rangle) \notag \\
&=\mathcal E \langle \phi_0,\phi_1 \rangle ( \abs{e_k(0)}+\abs{e_k(1)} ) \norm{e_k} \notag \\
&\;\;\;\; + \mathcal E \eta_k \langle \phi_0,\phi_1 \rangle (\abs{e_k(0)}+\abs{e_k(1)}) 
\end{align}
since $\langle \phi_0,\phi_1 \rangle=\langle \phi_{\mathcal N},\phi_{\mathcal N-1} \rangle$ by symmetry.

Thus, thanks to \eqref{e:ineg1}, \eqref{e:xx1}, \eqref{e:xx2}, \eqref{e:xx3}, \eqref{e:xx4} and \eqref{e:xx5}, the left-hand side of \eqref{e:etoile2} is greater than:
\begin{multline}
\label{e:xxx1}
\left( \frac{1}{\delt} + C_k \right)\norm{e_k}^2 - 
\left( \frac{2\eta_k+\norm{e_{k-1}}+ \mathcal E \langle \phi_0,\phi_1 \rangle ( \abs{e_k(0)}+\abs{e_k(1)} )}{\delt} + \sigma_k + f_k  \right) \norm{e_k} 
\\ - \frac{\eta_k \norm{e_{k-1}}+\mathcal E \eta_k \langle \phi_0,\phi_1 \rangle (\abs{e_k(0)}+\abs{e_k(1)})}{\delt} - \partieneg{C_k} \eta_k^2 - \eta_k f_k 
\end{multline}
\emph{Upper bound for RHS. }
We have:
\[ c(e_k,e_k,e_k^z) = c(e_k,e_k,e_k) - e_k(0) c(e_k,e_k,\phi_0) - e_k(1) c(e_k,e_k,\phi_{\mathcal N}) \]
but:
\begin{align*}
c(e_k,e_k,e_k) &= -\frac{1}{2} \int_0^1 e_k^2 \frac{\partial e_k}{\partial x} \\
&= - \frac{1}{6} \int_0^1 \frac{\partial \left[ (e_k)^3 \right] }{\partial x} \\
&= - \frac{1}{6} \left( \left(e_k(1)\right)^3 - \left(e_k(0)\right)^3 \right) 
\end{align*}
So:
\begin{align}
\label{e:xxx3}
c(e_k,e_k,e_k^z) &=  - \frac{1}{6} \left( \left(e_k(1)\right)^3 - \left(e_k(0)\right)^3 \right) - \left( e_k(0) \int_0^1 e_k^2 \frac{\partial \phi_0}{\partial x} + e_k(1) \int_0^1 e_k^2 \frac{\partial \phi_{\mathcal N}}{\partial x} \right) \notag\\
&=   - \frac{1}{6} \left( \left(e_k(1)\right)^3 - \left(e_k(0)\right)^3 \right) 
+ \mathcal N \left( e_k(0) \int_0^{1/\mathcal N} e_k^2 - e_k(1) \int_{1-1/\mathcal N}^1 e_k^2 \right)
\end{align}
Since, for all $x \in \left[ 0; \frac{1}{\mathcal N} \right]$,
\[ e_k(x) = e_k(0)+\mathcal N \left( e_k\left(\frac{1}{\mathcal N}\right) - e_k(0) \right) x \]
we have, thanks to $\abs{ e_k\left( \frac{1}{\mathcal N} \right) } \leq \mathcal E \norm{e_k} $:
\begin{align*}
\abs{ \mathcal N \times e_k(0) \int_0^{1/\mathcal N} e_k^2 } &\leq \mathcal N \frac{\abs{e_k(0)}}{\mathcal N} \left( \abs{e_k\left(\frac{1}{\mathcal N}\right)} \abs{e_k(0)} + \frac{ e_k\left(\frac{1}{\mathcal N}\right)^2 + e_k(0)^2}{3} + \frac{ 2 \abs{ e_k\left(\frac{1}{\mathcal N}\right) e_k(0) }}{3} \right)\\
&\leq \abs{e_k(0)} \Bigl(\mathcal E \norm{e_k} \abs{e_k(0)}  + \frac{\mathcal E^2 \norm{e_k}^2}{3}  +\frac{e_k(0)^2}{3} + \frac{2\abs{e_k(0)}}{3}\mathcal E \norm{e_k}\Bigr) \\
&\leq \frac{\mathcal E^2 \abs{e_k(0)}}{3} \norm{e_k}^2 + \frac{5}{3} \abs{e_k(0)}^2 \mathcal E \norm{e_k} + \frac{ \abs{e_k(0)}^3}{3}
\end{align*}
As a similar computation can be worked out for $ \abs{ \mathcal N \times e_k(1) \int_{1-1/\mathcal N}^{1} e_k^2 } $, we have:
\[ - \left( e_k(0) \int_0^1 e_k^2 \frac{\partial \phi_0}{\partial x} + e_k(1) \int_0^1 e_k^2 \frac{\partial \phi_{\mathcal N}}{\partial x} \right) \leq
 \xi_k^{\mathcal A} \norm{e_k}^2 + \xi_k^{\mathcal B} \norm{e_k} + \xi_k^{\gamma} \]
We also have, thanks to \eqref{e:triangle}:
\[ | r_k(e_k^z) | \leq \norm{r_k}_0 \norm{e_k^z} \leq \norm{r_k}_0 \norm{e_k} + \norm{r_k}_0 \eta_k \]
where:
\[ \norm{r_k}_0 = \sup_{v\in X_0, \norm v=1} r_k(v) \]
Hence, the right-hand side of \eqref{e:etoile2} is less than:
\[ \xi_k^{\mathcal A} \norm{e_k}^2 + ( \norm{r_k}_0 + \xi_k^{\mathcal B} ) \norm{e_k} + \frac{1}{6}\abs{e_k(1)^3-e_k(0)^3} + \xi_k^\gamma + \norm{r_k}_0\eta_k \]

\emph{Conclusion. }
Now \eqref{e:etoile2} implies, thanks to \eqref{e:xxx1}, and \eqref{e:xxx3}:
\begin{equation}
\label{e:etoile3}
\mathcal{A}_k \norm{e_k}^2 - 
   \mathcal{B}_k \norm{e_k} -
	\gamma_k \leq 0 
\end{equation}
Viewing left-hand side of \eqref{e:etoile3} as a (convex, thanks to our hypothesis \eqref{e:hypCn}) quadratic function $Q$ of $\norm{e_k}$, whose discriminant is $\mathcal D_k$, equation \eqref{e:etoile3} implies that, if $\mathcal D_k\geq0$, $\norm{e_k}$ is smaller than the greatest real root of $Q$, that is:
\[ \norm{e_k} \leq \frac{\mathcal B_k + \sqrt {\mathcal D_k}}{2 \mathcal A_k} \]
If $\mathcal D_k<0$, then necessarily $\gamma_k<0$ (as $\mathcal A_k$ is positive). Hence \eqref{e:etoile3} implies:
\[ \mathcal A_k \norm{e_k}^2 - 
   \mathcal B_k \norm{e_k} \leq 0 \]
and so:
\[ \norm{e_k} \leq \frac{\mathcal B_k}{\mathcal A_k}. \qedhere \]
\end{proof} 

\begin{acknowledgement} We would like to thank the anonymous referees, whose careful reading and comments have helped to greatly improve this paper. This work has been partially supported by the French National Research Agency (ANR) through COSINUS program (project COSTA-BRAVA n° ANR-09-COSI-015). \end{acknowledgement}

\bibliographystyle{plain}
\bibliography{biblio}

\end{document}